\documentclass[reqno]{amsart}%svjour2}
\usepackage{graphicx}
\usepackage{mathptmx}

\usepackage{amsmath}
\usepackage{amssymb}
\usepackage{amsfonts}
\usepackage{amsxtra}
\usepackage{color}

\usepackage{srcltx}
\numberwithin{equation}{section}

\def\soft#1{\leavevmode\setbox0=\hbox{h}\dimen7=\ht0\advance
    \dimen7 by-1ex\relax\if t#1\relax\rlap{\raise.6\dimen7
    \hbox{\kern.3ex\char'47}}#1\relax\else\if T#1\relax
    \rlap{\raise.5\dimen7\hbox{\kern1.3ex\char'47}}#1\relax
    \else\if d#1\relax\rlap{\raise.5\dimen7\hbox{\kern.9ex
    \char'47}}#1\relax\else\if D#1\relax\rlap{\raise.5\dimen7
    \hbox{\kern1.4ex\char'47}}#1\relax\else\if l#1\relax
    \rlap{\raise.5\dimen7\hbox{\kern.4ex\char'47}}#1\relax
    \else\if L#1\relax\rlap{\raise.5\dimen7\hbox{\kern.7ex
    \char'47}}#1\relax\else\message{accent \string\soft
    \space #1 not defined!}#1\relax\fi\fi\fi\fi\fi\fi}
\def\ocirc#1{\ifmmode\setbox0=\hbox{$#1$}\dimen0=\ht0
    \advance\dimen0 by1pt\rlap{\hbox to\wd0{\hss\raise\dimen0
    \hbox{\hskip.2em$\scriptscriptstyle\circ$}\hss}}#1\else
    {\accent"17 #1}\fi}
\def\XXint#1#2#3{{\setbox0=\hbox{$#1{#2#3}{\int}$ }
\vcenter{\hbox{$#2#3$ }}\kern-.6\wd0}}

\newtheorem{Theorem}[equation]{Theorem}
\newtheorem{Lemma}[equation]{Lemma}
\newtheorem{Proposition}[equation]{Proposition}
\newtheorem{Corollary}[equation]{Corollary}
 \newtheorem{Definition}[equation]{Definition}
\newtheorem{Remark}[equation]{Remark}

\newcommand{\norm}[1]{\lVert #1\rVert}
\newcommand{\bignorm}[1]{\big\lVert #1\big\rVert}
\newcommand{\setR}{\mathbb{R}}
\newcommand{\setN}{\mathbb{N}}
\newcommand{\embedding}{\hookrightarrow}
\newcommand{\compactembedding}{\hookrightarrow\hookrightarrow}
\newcommand{\bu}{{\bf u}}

\newcommand{\bg}{{\bf g}}
\newcommand{\bv}{{\bf v}}

\newcommand{\bW}{{\bf W}}

\newcommand{\ff}{{\bf f}}
\newcommand{\R}{\mathcal{R}}
\newcommand{\eps}{\epsilon}
\newcommand{\F}{\mathcal{F}}
\newcommand{\M}{\mathcal{M}}
\newcommand{\T}{\mathcal{T}}
\newcommand{\bnu}{\boldsymbol{\nu}}
\newcommand{\bnue}{\boldsymbol{\nu}_{\eta(t)}}
\newcommand{\bphi}{\boldsymbol{\phi}}
\newcommand{\bpsi}{\boldsymbol{\psi}}

\newcommand{\pa}{\partial}
\newcommand{\iot}{\int_{\Omega_{\eta(t)}}}

\newcommand{\iorte}{\int_{\Omega_{\mathcal{R}_\epsilon\eta(t)}}}

\newcommand{\iortd}{\int_{\Omega_{\mathcal{R}\delta(t)}}}
\newcommand{\ioretd}{\int_{\Omega_{\mathcal{R}_\eps\delta(t)}}}
\newcommand{\iortdn}{\int_{\Omega_{\mathcal{R}\delta_n(t)}}}
\newcommand{\idot}{\int_{\pa\Omega_{\eta(t)}\setminus\Gamma}}
\newcommand{\im}{\int_M}

\newcommand{\onorm}[1]{\norm{#1}_{L^2(\Omega_{\eta(t)})}}
\newcommand{\mnorm}[1]{\norm{#1}_{L^2(M)}}
\newcommand{\oet}{\Omega_{\eta(t)}}

\newcommand{\beweis}{\noindent{\bf Proof:}\ }

\newcommand{\loc}{\text{loc}}
\renewcommand{\phi}{\varphi}
\renewcommand{\P}{\mathcal{P}}
\newcommand{\bfzero}{\boldsymbol{0}}

\DeclareMathOperator{\dv}{div}

\DeclareMathOperator{\tr}{tr_\eta}
\DeclareMathOperator{\trt}{tr_{\tilde\eta}}

\DeclareMathOperator{\trnormal}{tr^n_\eta}
\DeclareMathOperator{\trnormaln}{tr^n_{\eta_0}}
\DeclareMathOperator{\trnormale}{tr^n_{\R_\epsilon\eta_0}}

\DeclareMathOperator{\trre}{tr_{\mathcal{R}_\epsilon\eta}}

\DeclareMathOperator{\trrd}{tr_{\mathcal{R}\delta}}
\DeclareMathOperator{\trred}{tr_{\mathcal{R}_\epsilon\delta}}
\DeclareMathOperator{\trrdt}{tr_{\mathcal{R}\delta(t)}}

\DeclareMathOperator{\grad}{grad}
\DeclareMathOperator{\id}{id}

\DeclareMathOperator{\spann}{span}
\DeclareMathOperator{\inn}{int}
\DeclareMathOperator{\supp}{supp}
\begin{document}

\title[Interaction of a generalized Newtonian fluid with a Koiter shell]{Global
  weak solutions for an incompressible, generalized Newtonian fluid interacting
  with a linearly elastic Koiter shell}
%[Interaction of a Newtonian fluid with a Koiter shell]
%[Interaction of a newtonian fluid with a Koiter shell]

%\titlerunning{Interaction of a Newtonian fluid with a Koiter shell}        % if too long for running head

\author{Daniel Lengeler}

%\authorrunning{} % if too long for running head

\address{Daniel Lengeler,
  Fakult{\"a}t f{\"u}r Mathematik\\
  Universit{\"a}t Regensburg\\
  Universit{\"a}tsstr.~31, 93040 Regensburg, Germany \\
} 
\email{daniel.lengeler@mathematik.uni-regensburg.de}

%\date{\notused}
% The correct date

\begin{abstract}
  In this paper we analyze the interaction of an incompressible, generalized 
  Newtonian fluid with a linearly elastic Koiter shell whose motion is
  restricted to transverse displacements. The middle surface of the
  shell constitutes the mathematical boundary of the three-dimensional
  fluid domain. We show that weak solutions exist as long as the
  magnitude of the displacement stays below some (possibly large)
  bound which is determined by the geometry of the undeformed shell.
% Insert your abstract here. Include up to five keywords.
\end{abstract}

\maketitle

\section{Introduction}
\label{intro}
Fluid-solid interaction problems involving moving interfaces have been
studied intensively during the last two decades. The interaction with
elastic solids has proven to be particularly difficult, due to
apparent regularity incompatibilities between the parabolic fluid
phase and the hyperbolic or dispersive solid phase, see, e.g., \cite{b37,b38,b40,b41,b60,b42,b27,b28,b61,b101} and the references therein. In \cite{b27,b28} the global-in-time existence of weak solutions for the interaction of an incompressible, Newtonian fluid with a Kirchhoff-Love plate is shown. In \cite{b101} we generalized this result to the case of a linearly elastic Koiter shell. The aim of the present paper is to extend the result in \cite{b101} to generalized Newtonian fluids, i.e., to fluids with a shear-dependent viscosity. A common model for the viscous (extra) stress tensor $S$ of such fluids is given by
\[S=\mu_0(\delta+|D|)^{p-2}D\]
for constants $\mu_0>0$, $\delta\ge 0$, $1<p<\infty$. Here, $D$ is the shear rate tensor. The mathematical analysis of such fluids in fixed spatial domains was initiated by Ladyzhenskaya \cite{b102,b103,b104} and Lions \cite{b105} in the late sixties. For $p\ge 11/5$ (in three space dimensions) the global existence of weak solutions follows from a combination of monotone operator theory and a compactness argument which is quite standard today. In \cite{b106} the Lipschitz trunction technique was used for the first time to study the existence of stationary weak solutions in the case of smaller exponents. This technique was improved in \cite{b107} and transfered to the nonstationary case in \cite{b108}. In the latter paper the existence of global weak solutions is shown for arbitary $p$ strictly greater than the natural bound $6/5$. This result is based on a parabolic Lipschitz truncation and a deep understanding of the pressure. However, in \cite{b109} the existence proof was considerably simplified by the introduction of solenoidal parabolic Lipschitz truncations which are considerably more flexible. We shall employ these in the present paper. It seems that \cite{b110} is the only analytical result so far dealing with the interaction of generalized Newtonian fluids with elastic solids. In this paper the existence of global weak solutions for shear-thickening fluids, i.e., $p\ge 2$, is shown under the assumptions of cylindrical symmetry, resulting in a two-dimensional problem, and a very strong mathematical damping of the elastic solid.

In the present paper we extend \cite{b101} to generalized Newtonian fluids. In doing so we have to deal with three new substantial difficulties. The first one is the well-known problem of identifying the limit of the extra stress tensor. Here, we have to apply the techniques developed in \cite{b109}. The second difficulty is due to the fact that the proof of relative $L^2$-compactness of bounded sequences of weak solutions developed in \cite{b101} needs substantial modification if $p$ is not larger than $3/2$. Finally, due to the additional nonlinearity in the system we cannot proceed as in \cite{b101} and apply the Kakutani-Glicksberg-Fan theorem. Instead, we have to construct an approximate decoupled system that is uniquely solvable on the one hand and that gives rise to an approximate coupled system accessible to the Lipschitz-trunction technique on the other hand. In order to deal with the approximate system we have to transfer  monotone operator theory techniques to the present ``non-cylindrical'' situation.

The present paper is partly based on the author's Ph.D. thesis \cite{b62}. It is organized as follows. In Subsection 1.1 we introduce Koiter's energy for elastic shells, in Subsection 1.2 we introduce the coupled fluid-shell system, and in Subsection 1.3 we derive formal a-priori estimates for this system. In Section 2 we give some results concerning domains with non-Lipschitz boundaries. Then, in Section 3 we state the main result of the paper. The rest of the section is devoted to the proof of this result. In Subsection 3.1 we give the proof of compactness of sequences of weak solutions. Subsequently, in Subsection 3.2 we analyse a decoupled variant of our original system, while in Subsection 3.3 we apply a fixed-point argument to this decoupled system. In Subsection 3.4 we conclude the proof by letting the regularisation parameter, which we introduced earlier, tend to zero. Finally, some further results and technical computations can be found in the appendix.

We write $W^{s,p}$ for the Sobolev-Slobodetskii scale of function spaces and, in particular, $H^s$ for the $L^2$-scale $W^{s,2}$. Furthermore, $\nabla$ denotes the Levi-Civita connection of Euclidean space or of given surfaces, depending on the context, and $\Delta$ is the corresponding Laplacian. Finally, we denote by $d\Phi$ the differential of mappings $\Phi$ between subsets of Euclidean space or of given surfaces.
\medskip

\subsection{Koiter's energy}
%%%%%%%%%%%%%%%%%%%%%%%%%%%%%%%%%%%%%%%%%%%%%%%%%%%%%%%%%%%%%%%%%%%%%%%%%%%%%%%%%%%%%%%%%%%%%%%%%%%%%%%%%%%%%%%%%%%%%%%%%%%%%%%%%%
%Koiter
Throughout the paper, let $\Omega\subset\setR^3$ be a bounded,
non-empty domain of class $C^4$ with outer unit normal $\bnu$. We
denote by $g$ and $h$ the first and the second fundamental form of
$\pa\Omega$, induced by the ambient Euclidean space, and by
$dA$ the surface measure of $\pa\Omega$. Furthermore, let
$\Gamma\subset\pa\Omega$ be a union of domains of class $C^{1,1}$
having non-trivial intersection with all connected components of
$\pa\Omega$. We set $M:=\pa\Omega\setminus\Gamma$; note that $M$
is compact. Let $\pa\Omega$ represent the middle surface of an
elastic shell of thickness $2\,\epsilon_0>0$ in its rest state where
$\epsilon_0$ is taken to be small compared to the reciprocal of the
principal curvatures. Furthermore, we assume that the elastic shell
consists of a homogeneous, isotropic material whose linear elastic
behavior may be characterized by the Lam{\'e} constants $\lambda$ and
$\mu$. We restrict the deformation of the middle surface to displacements along the unit normal field $\bnu$, and we assume the part $\Gamma$ of the middle surface to be fixed. Hence, we can describe the deformation by a scalar field $\eta:\ M
\rightarrow\setR$ vanishing at the boundary $\pa M$. We model the elastic energy of the deformation by \emph{Koiter's energy for linearly elastic shells and transverse
  displacements}
\begin{equation*}
 \begin{aligned}
  K(\eta)=K(\eta,\eta)=\frac{1}{2}\int_M \epsilon_0\,\langle
C,
\sigma(\eta)\otimes\sigma(\eta)\rangle + \frac{\epsilon_0^3}{3}\,\langle
C,
\xi(\eta)\otimes\xi(\eta)\rangle\ dA. 
 \end{aligned}
\end{equation*}
Here,
\begin{equation*}
  \begin{aligned}
    C_{\alpha\beta\gamma\delta}:=\frac{4\lambda\mu}{\lambda+2\mu} 
    \,g_{\alpha\beta}\,g_{\gamma\delta}
    + 2\mu\,(g_{\alpha\gamma}\,g_{\beta\delta} +
    g_{\alpha\delta}\,g_{\beta\gamma})
 \end{aligned}
\end{equation*}
is the elasticity tensor of the shell, and
\begin{equation*}
 \begin{aligned}
  \sigma(\eta)=-h\, \eta,\qquad  \xi(\eta)=\nabla^2 \eta - k\, \eta,
 \end{aligned}
\end{equation*}
are the linearized strain tensors, where $k_{\alpha\beta}:=h_\alpha^\sigma\, h_{\sigma\beta}$.
See \cite{b54}, \cite{b20}, \cite{b45}, \cite{b23} for Koiter's energy for nonlinearly elastic shells, and \cite{b23} for the derivation of the linearization; cf. also \cite{b101}. $K$ is a quadratic form in $\eta$ which is coercive on $H^2_0(M)$,
i.e., there exists a constant $c_0$ such that
\begin{equation}\label{eqn:kkoerziv}
 \begin{aligned}
K(\eta)\ge c_0\,\norm{\eta}_{H^2_0(M)}^2,
  \end{aligned}
\end{equation}
see the proof of Theorem 4.4-2 in \cite{b23}. Using integration by
parts and taking into account some facts from Riemannian geometry one
can show that the $L^2$-gradient of this energy has the form
\[\grad_{L^2}K(\eta)=\epsilon_0^3\,\frac{8\mu(\lambda+\mu)}{
  3(\lambda+2\mu)}\,\Delta^2\eta+B\eta\] where $B$ is a second order
differential operator which vanishes on flat parts of $M$, i.e., where
$h=0$. The details can be found in \cite{b62}. Thus, we obtain a generalization of the
linear Kirchhoff-Love plate equation for transverse displacements, cf.
for instance \cite{b44}.  By Hamilton's principle, the displacement
$\eta$ of the shell must be a stationary point of the action
functional
\[
\mathcal{A}(\eta)=\int_I \epsilon_0\rho_S\int_M
\frac{(\pa_t\eta(t,\cdot))^2}{2}\ dA-K(\eta(t,\cdot))\ dt
\]
where $I:=(0,T)$, $T>0$. Here we assume
that the mass density of $M$ may be described by a constant
$\epsilon_0\rho_S$. Hence, the integrand with respect to time is the
difference of the kinetic and the potential energy of the shell. The
corresponding Euler-Lagrange equation is
\begin{equation*}
  \begin{aligned}
    \epsilon_0\rho_S\,\pa^2_t\eta + \grad_{L^2}K(\eta)=0\text{ in }
    I\times M.
  \end{aligned}
\end{equation*}
\medskip

\subsection{Statement of the problem}\label{statement}

We denote by $\Omega_{\eta(t)}$, $t\in I$, the deformed domain
(cf.~\eqref{eq:def-om-eta}) and by \[\Omega_\eta^I:=\bigcup_{t\in I}\,
\{t\}\times \Omega_{\eta(t)}\] the deformed spacetime cylinder. Let us suppose that the variable domain $\Omega_\eta$ is filled by a
homogeneous, incompressible, generalized Newtonian fluid whose isothermal motion is governed by the system
\begin{equation}\label{eqn:fluid}
 \begin{aligned}
   \rho_F\big(\pa_t \bu + (\bu\cdot\nabla)\bu\big) - \dv \big(S(D\bu)
   - \pi\id\big) &= \rho_F\ff &&\mbox{ in }
   \Omega_{\eta}^I, \\
   \dv \bu &= 0 &&\mbox{ in }  \Omega_{\eta}^I,\\
   \bu(\,\cdot\,,\,\cdot\, + \eta\,\bnu) &= \pa_t\eta\,\bnu &&\mbox{ on } I\times M,\\
   \bu &= 0 &&\mbox{ on } I\times\Gamma.
\end{aligned}
\end{equation}
Here, $\bu$ is the velocity field, $\pi$ is the pressure field, $D\bu$ is the symmetric part of the gradient of $\bu$, $S$ is the extra stress tensor, $\id$ denotes the $3\times 3$ unit matrix, and $\ff$ is an external body force. We assume that $S$ possesses a \emph{$p$-structure}, i.e., for some $6/5<p<\infty$ and $\delta\ge 0$ we have
\begin{itemize}
 \item $S: M_{sym}\rightarrow M_{sym}$ continuous,
 \item Growth: $|S(D)|\le c_0(\delta+|D|)^{p-2}|D|$ for all $D\in M_{sym}$ and some $c_0>0$,
 \item Coercivity: $S(D):D\ge c_1(\delta+|D|)^{p-2}|D|^2$ for all $D\in
M_{sym}$ and some $c_1>0$,
\item Strict monotonicity: $(S(D)-S(E)):(D-E)>0$ for all $D,E\in M_{sym}, D\not=E$.
\end{itemize}
Here, $M_{sym}$ denotes the space of real, symmetric $3\times 3$ matrices. In the following, we divide equation \eqref{eqn:fluid}$_1$ by the constant fluid density $\rho_F$, denoting $S/\rho_F$ and $\pi/\rho_F$ again by $S$ and $\pi$. \eqref{eqn:fluid}$_{3,4}$ is the no-slip condition
in the case of a moving boundary, i.e., the velocity of the fluid at
the boundary equals the velocity of the boundary. The force exerted by
the fluid on the boundary is given by the evaluation of the stress
tensor at the deformed boundary in the direction of the inner normal
$-\bnu_{\eta(t)}$, i.e., by
\begin{equation}\label{eqn:kraft}
 \begin{aligned}
\rho_F\big(-S(D\bu(t,\cdot))\,\bnu_{\eta(t)} + \pi(t,\cdot)\,\bnu_{\eta(t)}\big).  
 \end{aligned}
\end{equation}
Thus, the equation for the displacement of the shell takes the form
\begin{equation}\label{eqn:shell}
\begin{aligned}
  \eps_0\rho_S\pa^2_{t} \eta + \grad_{L^2}K(\eta)  &= \eps_0\rho_S g + \rho_F{\bf F}\cdot\bnu
  &&\text{ in } I\times M, \\ 
  \eta=0,\ \nabla \eta &= \bfzero &&\text{ on } I\times\pa M
\end{aligned}
\end{equation}
where $g$ is a given body force and
\[{\bf F}(t,\cdot) := \big(-S(D\bu(t,\cdot))\,\bnu_{\eta(t)} +
\pi(t,\cdot)\,\bnu_{\eta(t)}\big)\circ\Phi_{\eta(t)}\ |\det
d\Phi_{\eta(t)}|\] with
$\Phi_{\eta(t)}:\pa\Omega\rightarrow\pa\Omega_{\eta(t)},\ q\mapsto
\eta(t,q)\,\bnu(q)$. In the following, we divide \eqref{eqn:shell}$_1$ by $\eps_0\rho_S$, denote $K/\eps_0\rho_S$ again by $K$, and assume, for the sake of a simple notation, that $\rho_F/\eps_0\rho_S=1$. Finally, we specify initial values
\begin{equation}\label{eqn:data}
 \begin{aligned}
   \eta(0,\cdot)=\eta_0,\ \pa_t\eta(0,\cdot)=\eta_1 \ \text{ in }M
   \quad\text{ and }\quad \bu(0,\cdot)=\bu_0 \ \text{ in
   }\Omega_{\eta_0}.
 \end{aligned}
\end{equation}
In the following, we will analyze the system \eqref{eqn:fluid}, \eqref{eqn:shell}, \eqref{eqn:data}. 
\medskip

\subsection{Formal a-priori estimates}\label{subsec:apriori}
Let us now formally derive energy estimates for this
parabolic-dispersive system. To this end, we multiply
\eqref{eqn:fluid}$_1$ by $\bu$, integrate the resulting identity over
$\Omega_{\eta(t)}$, and obtain after integrating by parts the stress
tensor\footnote{For the sake of a better readability we suppress the
  dependence of the unknown on the independent variables, e.g~we write
    $\bu=\bu(t,\cdot)$.}
\begin{align}
   &\iot \pa_t\bu\cdot\bu\ dx + \iot (\bu\cdot\nabla)\bu\cdot\bu\ dx
   \label{eqn:mult}
   \\
   &\quad= -\iot S(D\bu):D\bu\ dx+ \iot \ff\cdot\bu\ dx +\idot (S(D\bu)\,
   \bnue - \pi\,\bnue)\cdot \bu\ dA_{\eta(t)}.\notag
\end{align}
Here, $dA_{\eta(t)}$ denotes the surface measure of the deformed boundary $\pa\Omega_{\eta(t)}$. Taking into account that
\begin{equation}\label{eqn:wirbelid}
 \begin{aligned}
   \iot (\bu\cdot\nabla)\bu\cdot\bu\ dx=-\iot
   (\bu\cdot\nabla)\bu\cdot\bu\ dx + \int_{\pa\Omega_{\eta(t)}}
   \bu\cdot\bnu_{\eta(t)}|\bu|^2\ dA_{\eta(t)},
 \end{aligned}
\end{equation}
we may apply Reynold's transport theorem \ref{theorem:reynolds} to the first
two integrals in \eqref{eqn:mult} to obtain
\begin{equation}\label{eqn:energiefluid}
 \begin{aligned}
  \frac{1}{2}\frac{d}{dt} \iot |\bu|^2\ dx = &-\iot
S(D\bu):D\bu\ dx+ \iot \ff\cdot\bu\ dx\\
&+\idot
(S(D\bu)\, \bnue - \pi\,\bnue)\cdot \bu\ dA_{\eta(t)}.
 \end{aligned}
\end{equation}
% und verwenden
% wir die Koerzitit"at von $S$, so erhalten wir 
% \begin{equation}\label{ap:fluid}
%  \begin{aligned}
%   \frac{1}{2}\frac{d}{dt} &\iot |u|^2\ dx + c_1\iot
% |D\bu|^p\ dx - \iot\kappa_2\ dx\\ &-\idot
% (S(D\bu)\cdot \bnue + \pi\bnue)\cdot \bu\ dA_{\eta(t)} \le \iot \ff\cdot\bu\ dx.
%  \end{aligned}
% \end{equation}
Multiplying \eqref{eqn:shell}$_1$ by $\pa_t\eta$, integrating the
resulting identity over $M$, integrating by parts, and using the fact
that $(\grad_{L^2}K(\eta),\pa_t\eta)_{L^2}=2K(\eta,\pa_t\eta)$, we
obtain
\begin{equation}\label{eqn:energieshell}
 \begin{aligned}
   \frac{1}{2}\frac{d}{dt} \im |\pa_t\eta|^2\ dA + \frac{d}{dt}K(\eta)
   = \im g\, \pa_t\eta\ dA + \im {\bf F}\cdot\bnu\, \pa_t\eta\ dA.
 \end{aligned}
\end{equation}
% \begin{equation}\label{ap:shell}
%  \begin{aligned}
%   \frac{1}{2}\frac{d}{dt} \im |\pa_t\eta|^2\ dA + \frac{1}{2}\frac{d}{dt}\im &|\Delta\eta|^2\ dA +
% \im B\eta\ \pa_t\eta\ dA\\
%  &= \im g\ \pa_t\eta\ dA + \im {\bf F}\cdot\bnu \pa_t\eta\ dA.
%  \end{aligned}
% \end{equation}
% Offenbar gilt 
% \begin{equation*}
%  \begin{aligned}
%   \im B\eta\ \pa_t\eta\ dA \le c\big(\norm{\eta}_{H^2(M)}^2 + \norm{\pa_t \eta}_{L^2(M)}^2\big).
%  \end{aligned}
% \end{equation*}
Adding \eqref{eqn:energiefluid} and \eqref{eqn:energieshell}, taking
into account the definition of ${\bf F}$, \eqref{eqn:fluid}$_3$, and
applying a change of variables to the boundary integral, we obtain the
energy identity
\begin{equation}\label{eqn:energiesatz}
 \begin{aligned}
   \frac{1}{2}\frac{d}{dt} \iot |\bu|^2\ dx& + \frac{1}{2}\frac{d}{dt}
   \im |\pa_t\eta|^2\ dA +
   \frac{d}{dt}K(\eta)\\
   & = - \iot S(D\bu):D\bu\ dx + \iot \ff\cdot\bu\ dx + \im g\,
   \pa_t\eta\ dA.
 \end{aligned}
\end{equation}
In view of \eqref{eqn:kkoerziv} and the coercivity of $S$, an application of Gronwall's lemma
gives
% \begin{equation*}
%  \begin{aligned}
%    \frac{d}{dt}\onorm{\bu(t,\cdot)}^2 +
%    \norm{\nabla\bu(t,\cdot)}^2_{L^2(\Omega_{\eta(t)})}
%    +\frac{d}{dt}\mnorm{\pa_t\eta(t,\cdot)}^2 +
%    \frac{d}{dt}\norm{\eta(t,\cdot)}_{H^2(M)}^2
%    \\
%    \le c\,\big(\norm{\ff}_{L^2(\Omega_{\eta(t)})}^2 +
%    \norm{\bu}_{L^2(\Omega_{\eta(t)})}^2 + \norm{g}_{L^2(M)}^2+
%    \norm{\pa_t \eta}_{L^2(M)}^2\big).
%  \end{aligned}
% \end{equation*}
\begin{align}
  &\onorm{\bu(t,\cdot)}^2 +
  \int_0^t\norm{D\bu(s,\cdot)}^p_{L^p(\Omega_{\eta(s)})}\ ds +
  \mnorm{\pa_t\eta(t,\cdot)}^2 + \norm{\eta(t,\cdot)}_{H^2(M)}^2.\notag \\
  &\hspace{1cm}\le c\, e^t\,\Big(\norm{\bu_0}_{L^2(\Omega_{\eta_0})}^2 +
  \mnorm{\eta_1}^2 +
  \norm{\eta_0}_{H^2(M)}^2\label{ab:apriori} \\
&\hspace{4cm} + \int_0^t
  \norm{\ff(s,\cdot)}_{L^2(\Omega_{\eta(s)})}^2 + \mnorm{g(s,\cdot)}^2 ds\Big).\notag
 \end{align}
Hence, we have
\begin{equation*}
 \begin{aligned}
 \norm{\eta}_{W^{1,\infty}(I,L^2(M))\cap L^\infty(I,H_0^2(M))} +
\norm{\bu}_{L^\infty(I,L^2(\Omega_{\eta(t)}))}+\norm{D\bu}_{L^p(\Omega_\eta^I)}\le c(T,\text{data}).
 \end{aligned}
\end{equation*}
We shall construct weak solutions in this regularity class. In view of
the embedding $H^2(\pa\Omega)\embedding C^{0,\theta}(\pa\Omega)$ for
$\theta<1$, this implies that the boundary of our variable domain will
be the graph of a H{\"o}lder continuous function which, in general, is not
Lipschitz continuous. Since, in general, Korn's inequality is false in non-Lipschitz domains, c.f. \cite{b26}, we cannot expect an estimate of $\bu$ in $L^p(I,W^{1,p}(\Omega_{\eta(t)}))$. 
In the next section, we collect some facts about a class of non-Lipschitz domains.
\medskip

\section{Variable domains}
\label{sec:1}

%%%%%%%%%%%%%%%%%%%%%%%%%%%%%%%%%%%%%%%%%%%%%%%%%%%%%%%%%%%%%%%%%%%%%%%%%%%%%%%%%%%%%%%%%%%%%%%%%%%%%%%%%%%%%%%%%%%%%%%%%%%%%%%%%%
%Vargebiete
We denote by $S_{\alpha}$, $\alpha>0$, the open set of points in
$\setR^3$ whose distance from $\pa\Omega$ is less than $\alpha$. It's a
well known fact from elementary differential geometry, see for instance \cite{b1}, that there
exists a maximal $\kappa>0$ such that the mapping
\begin{equation*}
 \begin{aligned}
\Lambda: \partial\Omega\times (-\kappa,\kappa)\rightarrow
S_{\kappa},\
(q,s)\mapsto q + s\,\bnu(q)  
 \end{aligned}
\end{equation*}
is a $C^3$-diffeomorphism. For the inverse $\Lambda^{-1}$ we shall
write $x\mapsto(q(x),s(x))$. Note that $\kappa$ is not necessarily small; if $\Omega$ is the ball of radius $R$, then $\kappa=R$. Let $B_\alpha:=\Omega\cup S_{\alpha}$ for $0<\alpha<\kappa$. The
mapping $\Lambda(\,\cdot\, ,\alpha):\pa\Omega\rightarrow\pa B_\alpha$
is a $C^3$-diffeomorphism as well.  Hence, $B_\alpha$ is a bounded
domain with $C^3$-boundary.\footnote{In fact, it's even $C^4$.} For a continuous function
$\eta:\pa\Omega\rightarrow (-\kappa,\kappa)$ we set
\begin{equation}
 \begin{aligned}
  \Omega_\eta:=\Omega\setminus S_\kappa\ \cup \{x\in S_\kappa\ |\
  s(x)<\eta(q(x))\}.
 \end{aligned}\label{eq:def-om-eta}
\end{equation}
$\Omega_\eta$ is an open set. For $\eta\in
C^k(\pa\Omega)$, $k\in\{1,2,3\}$ we denote by $\bnu_\eta$ and
$dA_\eta$ the outer unit normal and the surface measure of
$\pa\Omega_\eta$, respectively. In \cite{b101} we showed that the mapping $\Psi_\eta:\overline\Omega\rightarrow\overline{\Omega_\eta}$, defined to be the identity in $\Omega\setminus S_\kappa$ and defined in $S_\kappa\cap\overline\Omega$ by
\begin{equation}\label{eqn:defpsi}
  \begin{aligned}
    &x\mapsto q(x)+\bnu(q(x))\big(s(x)+\eta(q(x))\,\beta(s(x)/\kappa)\big)
  \end{aligned}
\end{equation}
for a suitable function $\beta:\setR\rightarrow\setR$, is a homeomorphism, and even a $C^k$-diffeomorphism provided that $\eta\in C^k(\pa\Omega)$ with $k\in\{1,2,3\}$. Furthermore, we showed that the homeomorphism
\[\Phi_\eta:=\Psi_\eta|_{\pa\Omega}:\pa\Omega\rightarrow\pa\Omega_\eta,\
q\mapsto q+\eta(q)\,\bnu(q)\] is a $C^k$-diffeomorphism provided that $\eta\in C^k(\pa\Omega)$,
$k\in\{1,2,3\}$. Finally, we argued that $\Psi_\eta$ and $\Phi_\eta$ become singular as 
$\tau(\eta)\rightarrow\infty$ where
\begin{equation}\label{eqn:taueta}
  \begin{aligned}
    \tau(\eta):=\left\{\begin{array}{cl}
        (1-\norm{\eta}_{L^\infty(\pa\Omega)}/\kappa)^{-1} &
        \quad \text{ if }\norm{\eta}_{L^\infty(\pa\Omega)}<\kappa,\\
        \infty & \quad \text{ else.}\end{array}\right.
 \end{aligned}
\end{equation}
\smallskip
\begin{Remark}\label{bem:tdelta}
For $\eta\in C^2(\pa\Omega)$ with $\norm{\eta}_{L^\infty(\pa\Omega)}<\kappa$ and
$\bphi:\Omega\rightarrow\setR^3$ we denote by $\T_\eta\bphi$ the \emph{pushforward} of $(\det
d\Psi_\eta)^{-1} \bphi$ under $\Psi_\eta$, i.e., 
\[\T_\eta\bphi:=\big(d\Psi_\eta\, (\det
d\Psi_\eta)^{-1}\bphi\big)\circ\Psi^{-1}_\eta.\]
The mapping $\T_\eta$ with the inverse
$\T_\eta^{-1}\bphi:=\big(d\Psi_\eta^{-1}\, (\det
d\Psi_\eta)\,\bphi\big)\circ\Psi_\eta$
obviously defines isomorphisms between the Lebesgue and Sobolev spaces on $\Omega$ and
$\Omega_\eta$, respectively, as long as the order of differentiability is not larger than $1$. Moreover, the mapping preserves vanishing
boundary values. We saw in \cite{b101} that it also preserves the divergence-free constraint and hence defines isomorphisms between corresponding spaces of solenoidal functions
on $\Omega$ and $\Omega_\eta$, respectively.
\end{Remark}

A bi-Lipschitz mapping of domains induces isomorphisms of the
corresponding $L^p$ and $W^{1,p}$ spaces. For $\eta\in H^2(\pa\Omega)$
the mapping $\Psi_\eta$ is barely not bi-Lipschitz, due to the embedding
$H^2(\pa\Omega)\embedding C^{0,\theta}(\pa\Omega)$ for $\theta<1$.
Hence a small loss, made quantitative in the next lemma, will occur.
\begin{Lemma}\label{lemma:psi}
  Let $1<p\le\infty$ and $\eta\in H^2(\pa\Omega)$ with
  $\norm{\eta}_{L^\infty(\pa\Omega)}<\kappa$.  Then the linear mapping
  $v\mapsto v\circ\Psi_\eta$ is continuous from $L^p(\Omega_\eta)$ to
  $L^r(\Omega)$ and from $W^{1,p}(\Omega_\eta)$ to $W^{1,r}(\Omega)$
  for all $1\le r<p$. The analogous claim with $\Psi_\eta$ replaced by
  $\Psi_\eta^{-1}$ is true. The continuity constants depend only on
  $\Omega$, $p$, $r$, and a bound for $\norm{\eta}_{H^2(\pa\Omega)}$ and $\tau(\eta)$.
  % Konvergiert die Folge $(\eta_n)$ schwach gegen $\eta$ in
  % $H^2(\pa\Omega)$, so gilt f"ur $v\in W^{1,p}(\Omega_\eta)$
% \[v\circ\Psi_{\eta_n}\rightarrow v\circ\Psi\text{ in } W^{1,r}(\Omega).\]
\end{Lemma}
\proof See \cite{b101}.\qed \medskip

In the following we denote by $\ \cdot\ |_{\pa\Omega}$ the usual trace operator for Lipschitz domains. From the continuity properties of this trace operator and Lemma \ref{lemma:psi} we deduce the following assertion.
\begin{Corollary}\label{lemma:spur}
  Let $1<p\le\infty$ and $\eta\in H^2(\pa\Omega)$ with
  $\norm{\eta}_{L^\infty(\pa\Omega)}<\kappa$. Then the linear mapping
  $\tr: v\mapsto (v\circ\Psi_\eta)|_{\pa\Omega}$ is well defined and
  continuous from $W^{1,p}(\Omega_\eta)$ to $W^{1-1/r,r}(\pa\Omega)$
  for all $1<r<p$. The continuity constant depends only on $\Omega$, $r$, and a bound for
  $\norm{\eta}_{H^2(\pa\Omega)}$ and $\tau(\eta)$.
\end{Corollary}

From Lemma \ref{lemma:psi} and the Sobolev embeddings for regular domains we deduce Sobolev embeddings for our special domains.
\begin{Corollary}\label{lemma:sobolev}
  Let $1<p<3$ and $\eta\in H^2(\pa\Omega)$ with
  $\norm{\eta}_{L^\infty(\pa\Omega)}<\kappa$. Then\footnote{The symbol $\compactembedding$ indicates that the embedding is compact.}
  \[
  W^{1,p}(\Omega_\eta)\compactembedding L^s(\Omega_\eta)
  \] 
  for $1\le s < p^*=3p/(3-p)$. The embedding constant depends only on
  $\Omega$, $p$, $s$, and a bound for $\norm{\eta}_{H^2(\pa\Omega)}$ and $\tau(\eta)$.
\end{Corollary}

We denote by $H$ the mean curvature (with respect to the outer normal)
and by $G$ the Gauss curvature of $\pa\Omega$.

\begin{Proposition}\label{lemma:partInt}
  Let $1<p\le\infty$ and $\eta\in H^2(\pa\Omega)$ with
  $\norm{\eta}_{L^\infty(\pa\Omega)}<\kappa$.  Then, for $\bphi\in
  W^{1,p}(\Omega_\eta)$ with $\tr \bphi=b\,\bnu$, $b$ a scalar
  function, and $\psi\in C^1(\overline{\Omega_\eta})$ we have
  \[
  \int_{\Omega_\eta} \bphi\cdot\nabla\psi\ dx=-\int_{\Omega_\eta}
  \dv\bphi\ \psi\ dx+ \int_{\pa\Omega}b\, (1-2H\eta+G\,\eta^2)\, \tr
  \psi\ dA.
  \]
\end{Proposition}
\proof See \cite{b101}.\qed\medskip

We showed in \cite{b101} that the function $\gamma(\eta):=1-2H\eta+G\,\eta^2$ is positive as long as $|\eta|<\kappa$. Now, consider the space
\[E^p(\Omega_\eta):=\{\bphi\in
L^p(\Omega_\eta)\ |\ \dv \bphi\in
L^p(\Omega_\eta)\}\]
for $1\le p\le\infty$, endowed with the canonical norm. 
\begin{Proposition}\label{lemma:nspur}
  Let $1<p<\infty$ and $\eta\in H^2(\pa\Omega)$ with
  $\norm{\eta}_{L^\infty(\pa\Omega)}<\kappa$.  Then there exists a
  continuous, linear operator
\[\trnormal: E^p(\Omega_\eta)\rightarrow (W^{1,p'}(\pa\Omega))'\]
such that for $\bphi\in E^p(\Omega_\eta)$ and $\psi\in
C^1(\overline{\Omega_\eta})$
\[\int_{\Omega_\eta} \bphi\cdot\nabla\psi\ dx=-\int_{\Omega_\eta} \dv\bphi\
\psi\ dx + \langle
\trnormal\bphi,\tr\psi\rangle_{W^{1,p'}(\pa\Omega)}.\] The continuity
constant depends only on $\Omega$, $p$, and a bound for $\tau(\eta)$.
\end{Proposition}
\proof See \cite{b101}.\qed\medskip

\begin{Proposition}\label{lemma:FortVonRand}
  Let $1<p<\infty$, $\eta\in H^2(\pa\Omega)$ with
  $\norm{\eta}_{L^\infty(\pa\Omega)}<\kappa$, and $\alpha$ such that
  $\norm{\eta}_{L^\infty(\pa\Omega)}<\alpha<\kappa$. Then there exists a
  bounded, linear extension operator
  \[
  \F_\eta: \Big\{b\in W^{1,p}(\pa\Omega)\ \big|\
  \int_{\pa\Omega}b\,\gamma(\eta)\ 
  dA=0\Big\}\rightarrow W^{1,p}_{\dv}(B_\alpha),
  \]
  in particular $\tr\F_\eta b=b\,\bnu$. We also have
\[
  \F_\eta: \Big\{b\in L^{p}(\pa\Omega)\ \big|\ \int_{\pa\Omega}b\,
  \gamma(\eta)\ dA=0\Big\}\rightarrow \{\bphi\in L^p(B_\alpha)\ |\
  \dv\bphi=0\}
  \]
as a bounded, linear operator with $\trnormal\F_\eta b=b\,\gamma(\eta)$. The continuity constants depend only on $\Omega$, $p$, and a bound for $\norm{\eta}_{H^2(\pa\Omega)}$ and
  $\tau(\alpha)$.
\end{Proposition}
\proof See \cite{b101}.\qed\medskip

Of course, these extension operators are not optimal in the sense that they don't produce any regularity.

\begin{Proposition}\label{lemma:nullf}
  Let $6/5<p<\infty$ and $\eta\in H^2(\pa\Omega)$ with
  $\norm{\eta}_{L^\infty(\pa\Omega)}<\kappa$. Then extension by $0$
  defines a bounded, linear operator from $W^{1,p}(\Omega_\eta)$ to
  $H^{s}(\setR^3)$ for some $s>0$. The continuity constant depends only on
  $\Omega$, $p$, and a bound for $\norm{\eta}_{H^2(\pa\Omega)}$ and $\tau(\eta)$.
\end{Proposition}
\proof Let $6/5<r<p$. By standard embedding theorems we have $W^{1,r}(\Omega)\embedding H^{\tilde s}(\Omega)$ for some $\tilde s>0$. In order to prove the claim we can proceed exactly like in the proof of \cite[Proposition 2.28]{b101} once we showed that extension by $0$ defines a
bounded, linear operator from $H^{\tilde s}(\Omega)$ to $H^{s}(\setR^3)$ for some $0<s<\tilde s$.\footnote{In fact, it is possible to show that extension by $0$ is continuous from $H^s(\Omega)$ to $H^s(\setR^3)$ for all $0<s<1/2$. Hence, the application of H\"older's inequality below is not optimal.} To this end, it suffices to estimate the integral
  \begin{align}\notag
    \int_{\setR^3}\int_{\setR^3}& \frac{|v(x)-v(y)|^2}{|x-y|^{3+2s}}\
    dydx   = \int_{\Omega}\int_{\Omega}
    \frac{|v(x)-v(y)|^2}{|x-y|^{3+2s}}\ dydx +2
    \int_{\Omega}\int_{\setR^3\setminus\Omega}
    \frac{|v(x)|^2}{|x-y|^{3+2s}}\ dydx\\ 
    &\qquad=\int_{\Omega}\int_{\Omega} \frac{|v(x)-v(y)|^2}{|x-y|^{3+2s}}\
    dydx+ 2\int_{\Omega} |v(x)|^2 \int_{\setR^3\setminus\Omega}
    \frac{1}{|x-y|^{3+2s}}\ dydx \label{eqn:seuch}
 \end{align}
for $v\in H^{\tilde s}(\Omega)$. While the first term on the right-hand
side is dominated by $c\norm{v}_{H^{\tilde s}(\Omega)}$ for all $s\le \tilde s$, we can estimate
the interior integral of the second term by
\[
\int_{|z|>d(x)}\frac{1}{|z|^{3+2s}}\ dz=\frac{c(s)}{d(x)^{2s}}
\]
where $d(x)$ denotes the distance from $x$ to $\pa\Omega$. Again by standard embedding results, we have $H^{\tilde s}(\Omega)\embedding L^r(\Omega)$ for some $r>2$. An application of H\"older's inequality now shows that the second term on the right-hand side of \eqref{eqn:seuch} is dominated by \[c(s)\,\norm{v}^{2}_{L^r(\Omega)}\,\norm{d(\cdot)^{-2s}}_{L^{(r/2)'}(\Omega)}.\]
But the identity
\[
\int_{S_{\kappa/2}\cap\Omega}|d(x)|^{-2s(r/2)'}\
dx=\int_{\pa\Omega}\int_{-\kappa/2}^0|\det
d\Lambda|\, \alpha^{-2s(r/2)'} \ d\alpha dA,
\]
a consequence of a change of variables, proves that the last factor in this expression is finite for sufficiently small $s$.\qed\medskip

Let us now prove a suitable variant of Korn's inequality for non-Lipschitz domains.
\begin{Proposition}\label{prop:korn}
 Let $1<p<\infty$ and $\eta\in H^2(\pa\Omega)$ with
  $\norm{\eta}_{L^\infty(\pa\Omega)}<\kappa$. Then, for all $\bphi\in C^1(\overline{\Omega_\eta})$ and all $1\le r<p$, we have
\begin{equation}\label{eqn:korn}
 \norm{\bphi}_{W^{1,r}(\Omega_\eta)}\le c\, \big(\norm{D\bphi}_{L^{p}(\Omega_\eta)}+\norm{\bphi}_{L^{p}(\Omega_\eta)}\big).
\end{equation}
The constant $c$ depends only on $\Omega$, $p$, $r$, and a bound for $\norm{\eta}_{H^2(\pa\Omega)}$ and $\tau(\eta)$.
\end{Proposition}
\proof
For $1/\tilde r=1/r-1/p$ and $1- 1/(2\tilde r) <\beta<1$, we have
\[\norm{\nabla\bphi}_{L^{r}(\Omega_\eta)}\le \norm{\nabla\bphi\, d^{1-\beta}}_{L^{p}(\Omega_\eta)}\norm{d^{\beta-1}}_{L^{\tilde r}(\Omega_\eta)}.\]
Here, $d(x)$ denotes the distance from $x\in\Omega_\eta$ to $\pa\Omega_\eta$. Since $\Omega_\eta$ is a $\beta$-H\"older domain, by \cite[Theorem 3.1]{b100}, the term $\norm{\nabla\bphi\, d^{1-\beta}}_{L^{p}(\Omega_\eta)}$ is dominated by the right-hand side of \eqref{eqn:korn}. Hence, all we need to do is to bound the ${L^{\tilde r}(\Omega_\eta)}$-norm of $d^{\beta-1}$. To this end, we note that, since $\eta$ is $1/2$-H\"older continuous, for $|s|<\kappa$ and $q,\tilde q\in\pa\Omega$, we have
\begin{equation*}
 \begin{aligned}
|\eta(q)-s|\le |\eta(q)-\eta(\tilde q)| + |\eta(\tilde q)-s|&\le c\,|q-\tilde q|^{1/2} + |\eta(\tilde q)-s|\\
&\le c\,\big(|q-\tilde q| + |\eta(\tilde q)-s|\big)^{1/2}.
%\le c\, |q+s\bnu-\tilde q-\eta(\tilde q)\bnu|^{1/2}.
 \end{aligned}
\end{equation*}
For the second inequality, we used the fact that the geodesic distance on $\pa\Omega$ and the Euclidean distance in $\setR^3$ of $q$ and $\tilde q$ are comparable. We deduce that $d(q+s\bnu)\ge c\,|\eta(q)-s|^{2}$. Thus, using a change of variables, for $\norm{\eta}_{H^2(\pa\Omega)}<\alpha<\kappa$, we obtain
\begin{equation*}
 \begin{aligned}
\int_{S_\alpha\cap\Omega_\eta}d^{(\beta-1)\tilde r}\ dx &= \int_{\pa\Omega}\int_{-\alpha}^{\eta(q)}|\det
d\Lambda|\, d(q+s\bnu)^{(\beta-1)\tilde r}\ ds dA(q)\\
&\le c\,\int_{\pa\Omega}\int_{-\alpha}^{\alpha} |\eta(q)-s|^{2(\beta-1)\tilde r}\ dsdA(q). 
\end{aligned}
\end{equation*}
By the assumption on $\beta$ the last integral is bounded.
\qed
\medskip

The usual Bochner spaces are not the right objects to deal with
functions defined on time-dependent domains. For this reason we now
define an (obvious) substitute for these spaces. For $I:=(0,T)$,
$T>0$, and $\eta\in C(\bar I\times\pa\Omega)$ with
$\norm{\eta}_{L^\infty(I\times\pa\Omega)}<\kappa$ we set
$\Omega_\eta^I:=\bigcup_{t\in I}\, \{t\}\times \Omega_{\eta(t)}$.
Note that $\Omega_\eta^I$ is a domain in $\setR^4$. For $1\le
p,r\le\infty$ we set
\begin{equation*}
 \begin{aligned}
   L^p(I,L^r(\oet))&:=\{v\in L^1(\Omega_\eta^I)\ |\ v(t,\cdot)\in
   L^r(\Omega_{\eta(t)})\text{ for almost all $t$ and}\\
   &\hspace{.8cm}\norm{v(t,\cdot)}_{L^r(\Omega_{\eta(t)})}\in
   L^p(I)\},\\
   L^p(I,W^{1,r}(\Omega_{\eta(t)}))&:=\{v\in L^p(I,L^r(\oet))\ |\
   \nabla v\in
   L^p(I,L^r(\oet))\},\\
   L^p(I,W_{\dv}^{1,r}(\Omega_{\eta(t)}))&:=\{\bv\in
   L^p(I,W^{1,r}(\Omega_{\eta(t)}))\ |\
   \dv \bv=0\},\\
L^p(I,W_{\dv,s}^{1,r}(\Omega_{\eta(t)}))&:=\{\bv\in
   L^p(I,L^r(\Omega_{\eta(t)}))\ |\
   D\bu\in L^p(I,L^r(\Omega_{\eta(t)})),\, \dv \bv=0\},\\
   W^{1,p}(I,W^{1,r}(\Omega_{\eta(t)}))&:=\{\bv\in
   L^p(I,W^{1,r}(\Omega_{\eta(t)}))\ |\ \pa_t\bv\in
   L^p(I,W^{1,r}(\Omega_{\eta(t)}))\}.
 \end{aligned}
\end{equation*}
Here $\nabla$ and $\dv$ are acting with respect to the space variables. Furthermore, we set
\[
\Psi_\eta: \bar I\times \overline\Omega\rightarrow
\overline{\Omega_\eta^I},\ (t,x)\mapsto (t,\Psi_{\eta(t)}(x))
\]
and
\[
\Phi_\eta: \bar I\times \pa\Omega\rightarrow\bigcup_{t\in \bar I}\,
\{t\}\times \pa\Omega_{\eta(t)},\ (t,x)\mapsto (t,\Phi_{\eta(t)}(x)).
\] 
If $\eta\in L^\infty(I,H^2(\pa\Omega))$ we obtain
``instationary'' versions of the claims made so far by applying these
at (almost) every $t\in I$. For instance, from Corollary \ref{lemma:sobolev} we
deduce that
\[
L^2(I,H^{1}(\Omega_{\eta(t)}))\embedding L^2(I,L^s(\Omega_{\eta(t)}))
\]
for $1\le s<2^*$. Note that the construction given above does not
provide a substitute for Bochner spaces of functions with values in negative spaces.
% Diese wird "ublicherweise durch lokale Trivialisierungen induziert. Das
% sind Abbildungen, die das B"undel lokal auf das Kreuzprodukt einer offenen Teilmenge der
% Mannigfaltigkeit mit der sogenannten typischen Faser, einem festen Vektorraum, abbilden. 
% ergibt sich durch
% Anwenden des obigen Spuroperators zu jedem Zeitpunkt ein stetiger Operator $\tr$ von
% $L^p(I,W^{1,p}(\Omega_\eta(t)))$ nach $L^p(I,W^{1-1/r,r}(\pa\Omega))$ f"ur alle $r<p$.
% \marginpar{PR"UFEN!}
Furthermore, note that for all $1/2<\theta<1$ we have
\begin{align}
  \begin{aligned}
    &W^{1,\infty}(I,L^2(\pa\Omega))\cap
    L^\infty(I,H^2(\pa\Omega))\label{eqn:hoeldereinb}
    \\
    &\hspace{2cm} \embedding C^{0,1-\theta}(\bar
    I,H^{2\theta}(\pa\Omega))\embedding C^{0,1-\theta}(\bar I,
    C^{0,2\theta -1}(\pa\Omega)).
  \end{aligned}
\end{align}

\begin{Proposition}\label{lemma:FortVonRandZeit}
  Let $\eta\in W^{1,\infty}(I,L^2(\pa\Omega))\cap
  L^\infty(I,H^2(\pa\Omega))$ be given with \linebreak
  ${\norm{\eta}_{L^\infty(I\times\pa\Omega)}<\kappa}$ and $\alpha$ a
  real number such that
  $\norm{\eta}_{L^\infty(I\times\pa\Omega)}<\alpha<\kappa$. The
  application of the extension operators from Proposition
  \ref{lemma:FortVonRand} at (almost) all times defines a bounded,
  linear extension operator $\F_\eta$ from
  \[
  \Big\{b\in H^{1}(I,L^2(\pa\Omega))\cap L^2(I,H^2(\pa\Omega))\ |\
  \int_{\pa\Omega}b(t,\cdot)\, \gamma(\eta(t,\cdot))\ dA=0\text{ for
    all }t\in I\Big\}
  \] 
  to
  \[
  \{\bphi\in H^{1}(I,L^2(B_\alpha))\cap C(\bar I,H^1(B_\alpha))\ |\
  \dv\bphi=0\},
  \]
as well as a bounded,
  linear extension operator $\F_\eta$ from
  \[
  \Big\{b\in C(\bar I,L^2(\pa\Omega))\ |\ \int_{\pa\Omega}b(t,\cdot)\,
  \gamma(\eta(t,\cdot))\ dA=0\text{ for almost all }t\in I\Big\}
  \] 
  to
  \[
  \{\bphi\in C(\bar I,L^2(B_\alpha))\ |\ \dv\bphi=0\}.
  \]
  The continuity constants depend only on $\Omega$ and a bound for $\norm{\eta}_{
    W^{1,\infty}(I,L^2(\pa\Omega))\cap L^\infty(I,H^2(\Omega))}$ and $\tau(\alpha)$.
\end{Proposition}
\proof See \cite{b101}.\qed\smallskip

\begin{Remark}\label{bem:tdelta2}
For $\eta\in C^2(I\times\pa\Omega)$ with $\norm{\eta}_{L^\infty(I\times\pa\Omega)}<\kappa$ an application of $\T_{\eta(t)}$
for each $t\in I$ defines isomorphisms between appropriate function spaces on $I\times\Omega$ and $\Omega_\eta^I$, respectively, as
long as the order of differentiability is not larger than $1$.
\end{Remark}

\section{Main result}
\label{sec:2}

%%%%%%%%%%%%%%%%%%%%%%%%%%%%%%%%%%%%%%%%%%%%%%%%%%%%%%%%%%%%%%%%%%%%%%%%%%%%%%%%%%%%%%%%%%%%%%%%%%%%%%%%%%%%%%%%%%%%%%%%%%%%%%%%%%
%Existenz
For the rest of the paper we shall fix some $6/5<p<\infty$. We define 
\begin{equation*}
 Y^I:=W^{1,\infty}(I,L^2(M))\cap L^\infty(I,H_0^2(M)),
\end{equation*}
and for $\eta\in Y^I$ with $\norm{\eta}_{L^\infty(I\times M)}<\kappa$ we set
\begin{equation*}
X_{\eta,p}^I:=L^\infty(I,L^2(\Omega_{\eta(t)}))\cap L^p(I,W^{1,p}_{\dv,s}(\oet)).
\end{equation*}
Here and throughout the rest of the paper, we tacitly extend functions
defined in $M$ by $0$ to $\pa\Omega$. Note that, by Proposition \ref{prop:korn}, the space $X_{\eta,p}^I$ embeds into $L^r(I,W^{1,r}_{\dv}(\oet))$ for all $1\le r<p$.  We define the space of test functions $T_{\eta,p}^I$ to consist of  all couples
\[(b,\bphi)\in \big(H^1(I,L^2(M))\cap L^{\tilde p}(I,H^2_0(M))\big)\times
\big(H^1(I,L^2(\oet))\cap
L^{\tilde p}(I,W^{1,\tilde p}_{\dv}(\oet))\big)\]
such that $b(T,\cdot)=0$, $\bphi(T,\cdot)=0$,\footnote{We saw in \cite{b101} that it makes sense to evaluate $\bphi$ at a fixed point $t$ in time and that $\bphi(t,\cdot)\in L^2(\oet)$.}
% \footnote{Der Divergenzoperator wirkt wie
% immer nur auf die r"aumlichen Koordinaten.}
and $\bphi-\F_\eta b\in H_0$. Here, $H_0$ denotes the closure in $H^1(I,L^2(\oet))\cap
L^{\tilde p}(I,W^{1,\tilde p}_{\dv}(\oet))$ of the elements of this space that vanish at $t=T$ and whose supports are contained in $\Omega_{\eta}^{\bar I}$. From the last requirement we infer that $\tr\bphi=\tr \F_\eta b=b\,\bnu$. In particular, $\bphi$ vanishes on $\Gamma$. Furthermore, the finite exponent $\tilde p$ needs to be larger than $(5p/6)'$ and not smaller than $p$, so let us choose $\tilde p:=\max((5p/6)'+3,p)$. 

We call the data $(\ff,g,\bu_0,\eta_0,\eta_1)$ \emph{admissible} if $\ff\in L_{\text{loc}}^2([0,\infty)\times
\setR^3)$, $g\in L_{\text{loc}}^2([0,\infty) \times M)$, $\eta_0\in H^2_0(M)$ with
$\norm{\eta_0}_{L^\infty(M)}<\kappa$, $\eta_1\in L^2(M)$, and $\bu_0\in
L^2(\Omega_{\eta_0})$ with $\dv \bu_0=0$,
$\trnormaln\bu_0=\eta_1\,\gamma(\eta_0)$.

\begin{Definition} A couple $(\eta,\bu)$ is a weak solution of
  \eqref{eqn:fluid}, \eqref{eqn:shell}, and \eqref{eqn:data} for the
  admissible data $(\ff,g,\bu_0,\eta_0,\eta_1)$ in the intervall $I$
  if $\eta\in Y^I$ with $\norm{\eta}_{L^\infty(I\times M)}<\kappa$,
  $\eta(0,\cdot)=\eta_0$, $\bu\in X_{\eta,p}^I$ with $\tr \bu =
  \pa_t\eta\,\bnu$, and
  \begin{align}
    &- \int_I\iot \bu\cdot\pa_t\bphi\ dxdt +\int_I\iot\bu\otimes\bu:D\bphi\
    dxdt + \int_I\iot S(D\bu):D\bphi\ dxdt\notag\\
& -\int_I\im\pa_t\eta\,
    \pa_tb\ dAdt + 2\int_I K(\eta,b)\ dt\label{eqn:schwach}
    \\
    &\hspace{0.5cm}=\int_I\iot \ff\cdot\bphi\ dxdt + \int_I\im g\, b\
    dAdt+\int_{\Omega_{\eta_0}}\bu_0\cdot\bphi(0,\cdot)\ dx +
    \im\eta_1\, b(0,\cdot)\ dA \notag
  \end{align}
  for all test functions $(b,\bphi)\in T_{\eta,p}^I$.
\end{Definition}

Like in \cite{b101} the weak formulation \eqref{eqn:schwach} arises formally by
multiplication of \eqref{eqn:fluid} with a test function $\bphi$,
integration over space and time, integration by parts, and taking into
account \eqref{eqn:shell}. Here, the boundary integrals resulting from integrating by parts the time-derivative of $\bu$ and the convective term cancel. By Corollary \ref{lemma:sobolev} and interpolation (with a weight of $\theta=2/5$ on the bound for the kinetic energy), we have $\bu\in L^r(\Omega_\eta^I)$ for all $1\le r<10p/6$. Hence, in view of the assumption on $p$, the second term in \eqref{eqn:schwach} is well-defined and finite.

\begin{Theorem}\label{theorem:hs}
  For arbitrary admissible data $(\ff,g,\bu_0,\eta_0,\eta_1)$ there
  exist a time $T^*\in (0,\infty]$ and a couple $(\eta,\bu)$ such that
  for all $T<T^*$ $(\eta,\bu)$ is a weak solution of
  \eqref{eqn:fluid}, \eqref{eqn:shell}, and \eqref{eqn:data} in the
  intervall $I=(0,T)$. Furthermore, we have
  \begin{align}
      &\norm{\eta}_{Y^I}^2 + \norm{\bu}_{L^\infty(I,L^2(\Omega_{\eta(t)}))}^2 + \norm{D\bu}_{L^p(\Omega_{\eta}^I)}^p\label{ab:hs}\\\notag
&\le c\,e^T\Big(\norm{\bu_0}_{L^2(\Omega_{\eta_0})}^2 + \mnorm{\eta_1}^2 +
  \norm{\eta_0}_{H^2(M)}^2 + \int_0^T
  \norm{\ff(s,\cdot)}_{L^2(\Omega_{\eta(s)})}^2 + \mnorm{g(s,\cdot)}^2 ds\Big).
  \end{align}
Either $T^*=\infty$ or $\lim_{t\rightarrow T^*}\norm{\eta(t,\cdot)}_{L^\infty(M)}=\kappa$.
\end{Theorem}

In the following we will denote the right-hand side of \eqref{ab:hs}
as a function of $T$, $\Omega_\eta^I$, and the data by
$c_0(T,\Omega_\eta^I,\ff,g,\bu_0,\eta_0,\eta_1)$.

\subsection{Compactness}
\label{subsec:21}
Similarly to \cite{b101} we can show strong $L^2$-compactness of the shell and the fluid velocities for bounded sequences of weak solutions. However, for the compactness of the shell velocities we need to assume that $p>3/2$. The reason is that we need the shell velocities to be uniformly bounded in a spatial regularity class that embeds compactly into $L^2(M)$. By taking the trace of the fluid velocities, we obtain the boundedness of the shell velocities in $L^p(I,W^{1-1/r,r}(M))$ for all $1\le r<p$. But $W^{1-1/r,r}(M)$ embeds compactly into $L^2(M)$ if and only if $r>3/2$. While the weak formulation \eqref{eqn:schwach} of our original system is linear in the shell velocity (and compactness of the shell velocities is therefore not needed), this is not the case in our regularized system. On the other hand, since the extra stress tensor of our regularized system will possess a $p$-structure for some large $p$ (partly in order to make the problem accessible to monotone operator theory), in the end, we can deal with arbitrary $p>6/5$.

\begin{Proposition}\label{lemma:komp}
Let $(\ff,g,\bu_0^n,\eta_0^n,\eta_1^n)$ a sequence of admissible data with
\begin{equation}\label{eqn:abeschraenkt}
 \begin{aligned}
\sup_n\big(\tau(\eta_0^n)+\norm{\eta_0^n}_{H^2_0(M)}+\norm{\eta_1^n}_{L^2(M)}+\norm{
\bu^n_0}_{
L^2(\Omega_ {\eta_0^n})}\big)<\infty.
 \end{aligned}
\end{equation} 
Furthermore, let $(\eta_n,\bu_n)$ be a sequence of weak solutions of \eqref{eqn:fluid},
\eqref{eqn:shell}, and \eqref{eqn:data} for the above data in the intervall $I=(0,T)$ such that
% \,\footnote{F"ur zeitabh"angige
% Funktionen $\eta$ ersetzen wir  in
% der Definition von $\tau(\eta)$ die Norm $\norm{\eta}_{L^\infty(\pa\Omega)}$ durch
% $\norm{\eta}_{L^\infty(I\times \pa\Omega)}$.}
\begin{equation}\label{ab:beschraenkt}
 \begin{aligned}
\sup_n\big(\tau(\eta_n) + \norm{\eta_n}_{Y^I} +
\norm{\bu_n}_{X_{\eta_n,p}^I}\big)<\infty.
 \end{aligned}
\end{equation}
Then the sequence $(\bu_n)$ is relatively compact in $L^2(I\times
\setR^3)$.\footnote{Here and throughout the rest of the paper, if not explicitly stated otherwise, we (tacitly) extend functions defined in a domain of $\setR^3$ by $0$ to the whole space.} If $p>3/2$, the sequence $(\pa_t\eta_n)$ is relatively compact in $L^2(I\times M)$.
% \footnote{Wie immer setzen wir die Felder $\bu_n$ durch
% $\boldsymbol{0}$ auf $I\times\setR^3$ fort.}
\end{Proposition}
\proof We infer from \eqref{ab:beschraenkt} that for a subsequence\footnote{When passing over to a subsequence we will
tacitly always do so with respect to all involved sequences and use again the subscript $n$.} we have
 \begin{alignat}{2}\notag
   \eta_n&\rightarrow\eta &&\quad\text{ weakly$^*$ in
   }L^\infty(I,H^2_0(M))\text{ and uniformly},\\ \label{eqn:schw}
   \pa_t\eta_n&\rightarrow\pa_t\eta &&\quad\text{ weakly$^*$ in }L^\infty(I,L^2(M)),\\\notag
   \bu_n&\rightarrow\bu &&\quad\text{ weakly$^*$ in } L^\infty(I,L^2(\setR^3)),\\\notag
   \nabla\bu_n&\rightarrow\nabla\bu &&\quad\text{ weakly in } L^2(I\times \setR^3).
\end{alignat}
Here, we extend the functions $\nabla\bu_n$ and $\nabla\bu$, which a-priori are defined
only in $\Omega_{\eta_n}^I$ and $\Omega_{\eta}^I$, respectively, by $\bfzero$ to $I\times \setR^3$. Let us deal with the case $p>3/2$ first. The proof of this case is a rather simple modification of the proof of \cite[Proposition 3.8]{b101}. Therefore, we only give a sketch. We saw in the proof of \cite[Proposition 3.8]{b101} that it is enough to show that
 \begin{align}\notag
 \int_I\int_{\Omega_{\eta_n(t)}}\bu_n\cdot \F_{\eta_n}\pa_t\eta_n\ dxdt
&+\int_I\im|\pa_t\eta_n|^2\ dAdt\\\label{eqn:l2konv1}
\rightarrow \int_I\int_{\Omega_{\eta(t)}}\bu\cdot \F_{\eta}\pa_t\eta\ dxdt 
&+ \int_I\im|\pa_t\eta|^2\ dAdt,\\\notag
\int_I\int_{\Omega_{\eta_n(t)}} \bu_n\cdot(\bu_n -
\F_{\eta_n}\pa_t\eta_n)\ dxdt &\rightarrow\int_I \int_{\Omega_{\eta(t)}} \bu\cdot(\bu -
\F_{\eta}\pa_t\eta)\ dxdt.
 \end{align}
Here, we assume that the number $\alpha$ in the definition of $\F$, see
Proposition \ref{lemma:FortVonRandZeit}, satisfies the inequality
$\sup_n\norm{\eta_n}_{L^\infty(I\times M)}<\alpha<\kappa$. Let us start with the demonstration of \eqref{eqn:l2konv1}$_1$. For $b\in H^2_0(M)$ we employ the special test functions $(\M_{\eta_n}b,\F_{\eta_n}\M_{\eta_n}b)$, see Lemma \ref{lemma:mittelwert} for the definition of the operators $\M_{\eta_n}$. From this lemma, Proposition \ref{lemma:FortVonRandZeit}, Proposition \ref{lemma:FortVonRand}, and \eqref{ab:beschraenkt} we deduce the estimate
\begin{equation*}
 \begin{aligned}
   \norm{\M_{\eta_n}b}_{H^1(I,L^2(M))\cap
     L^{\tilde p}(I,H^2_0(M))}+\norm{\F_{\eta_n}\M_{\eta_n}b}_{H^1(I,L^2(B_\alpha))\cap
     C(\bar I,H^1(B_\alpha))\cap L^{\tilde p}(I,W^{1,\tilde p}(B_\alpha))}\\
   \le c\,\norm{b}_{H^2_0(M)}.
 \end{aligned}
\end{equation*}
As in the proof of \cite[Proposition 3.8]{b101} we use equation \eqref{eqn:schwach} to show that the functions 
\begin{align*}
c_{b,n}(t)&:=\int_{\Omega_{\eta_n(t)}}\bu_n(t,\cdot)\cdot(\F_{\eta_n}\M_{\eta_n}b)(t,\cdot)\ dx+\int_M\pa_t\eta_n(t,\cdot)\,(\M_{\eta_n}b)(t,\cdot)\ dA
\end{align*}
are bounded in $C^{0,\beta}(\bar I)$ for some $\beta\in (0,1)$ independently of $\norm{b}_{H^2_0(M)}\le 1$. Here, the convective term has to be estimated in the form
\begin{align*}
\bignorm{\int_{\Omega_{\eta_n(t)}}\bu_n\otimes\bu_n:D\F_{\eta_n}\M_{\eta_n}b\ dx}_{L^{1/(1-\beta)}(I)}\le \norm{\bu_n}_{L^{2\tilde p'}(\Omega_{\eta_n}^I)}^2\,
\norm{D\F_{\eta_n}\M_{\eta_n}b}_{L^{\tilde p}(\Omega_{\eta_n}^I)}\\
\le c\norm{\bu_n}_{L^\infty(I,L^2(\Omega_{\eta_n(t)}))}^{4/5}\norm{\bu_n}_{L^p(I,W^{1,p}(\Omega_{\eta_n(t)}))}^{6/5}\norm{D\F_{\eta_n}\M_{\eta_n}b}_{L^{\tilde p}(\Omega_{\eta_n}^I)}.
\end{align*}
Note that $2\tilde p'<10p/6$. From this fact and \eqref{ab:beschraenkt} we deduce as before by the Arzela-Ascoli argument that the functions
\begin{equation}\label{eqn:hn}
h_n(t):=\sup_{\norm{b}_{H^2_0(M)}\le
1} \big(c_{b,n}(t)-c_{b}(t)\big), 
\end{equation}
where $c_{b}$ is defined as $c_{b,n}$ with $(\eta_n,\bu_n)$ replaced by $(\eta,\bu)$, converge to zero in $C(\bar I)$. By \cite[Lemma A.13]{b101}, for the functions
\[g_n(t):=\sup_{\norm{b}_{L^2(M)}\le
1} \big(c_{b,n}(t)-c_{b}(t)\big)\]
and all $3/2< r<p$ we have
\begin{equation}\label{eqn:gest}
\int_I g_n(t)\ dt \le \epsilon\,
c\,\big(\norm{\bu_n}_{L^p(I,W^{1,r}(\Omega_{\eta_n(t)}))}+\norm{\bu}_{L^p(I,W^{1,r}(\Omega_{\eta(t)}))}
\big) + c(\epsilon)\int_I h_n(t)\ dt, 
\end{equation}
proving that $(g_n)$ tends to zero in $L^1(I)$. As in the proof of \cite[Proposition 3.8]{b101} we can infer \eqref{eqn:l2konv1}$_1$. Let us proceed with the proof of \eqref{eqn:l2konv1}$_2$. We fix a sufficiently small $\sigma>0$ and $\delta_\sigma\in C^4(\bar I\times \pa\Omega)$ such that $\norm{\delta_\sigma-\eta}_{L^\infty(I\times
\pa\Omega)}<\sigma$ and $\delta_\sigma<\eta$ in $\bar I\times\pa\Omega$. For $\bphi\in
H(\Omega)$ and $t\in\bar I$ we set
\[c_{\bphi,n}^\sigma(t):=\int_{\Omega_{\eta_n(t)}}\bu_n(t,\cdot)\cdot \T_{\delta_\sigma(t)}\bphi\
dx,\]
see Remark \ref{bem:tdelta} for the definition of $\T_{\delta_\sigma(t)}$, and we define the functions $c_{\bphi}^\sigma$ analogously. From Remark \ref{bem:tdelta} and Remark \ref{bem:tdelta2} we deduce that
\[\norm{\T_{\delta_\sigma}\bphi}_{H^1(I,L^2(B_\alpha))\cap C(\bar I,W^{1,\tilde p}(B_\alpha))}\le
c\,\norm{\bphi}_{W^{1,\tilde p}_{0,\dv}(\Omega)}.\]
As before, using equation \eqref{eqn:schwach}, we infer that the functions
\[h_n^\sigma(t):=\sup_{\norm{\bphi}_{W^{1,\tilde p}_{0,\dv}(\Omega)}\le
1} \big(c_{\bphi,n}^\sigma(t)-c_{\bphi}^\sigma(t)\big)\]
are bounded in some H\"older space, independently of $\norm{b}_{H^2_0(M)}$. Again by the Arzela-Ascoli argument, we obtain that $(h_n^\sigma)$ tends to zero in $C(\bar I)$, and, by an application of \cite[Lemma A.13]{b101}, that the functions
\[g_n^\sigma(t):=\sup_{\norm{\bphi}_{H(\Omega)}\le
1} \big(c_{\bphi,n}^\sigma(t)-c_{\bphi}^\sigma(t)\big)\]
converge to zero in $L^1(I)$. Finally, \cite[Lemma A.16]{b101} yields the existence of functions $\bpsi_{t,n}$ as in the proof of \cite[Proposition 3.8]{b101} satisfying the estimate
\begin{equation*}
  \norm{\bu_n(t,\cdot)-(\F_{\eta_n}\pa_t\eta_n)(t,\cdot) -
\bpsi_{t,n}}_{(H^{s}(\setR^3))'} <
\epsilon,
\end{equation*}
for arbitray, but fixed $s>0$. On the other hand, by Lemma \ref{lemma:nullf} the functions $\bu_n$ and $\bu$, extended by $\bfzero$ to $I\times\setR^3$, are uniformly bounded in $L^p(I,H^s(\setR^3))$ for sufficiently small $s$. Thus, we can infer \eqref{eqn:l2konv1}$_2$ as in the proof of \cite[Proposition 3.8]{b101}.

Now, let us consider the case $6/5<p\le 3/2$. In view of \eqref{eqn:schw}$_3$ we have
\[\limsup_{n\rightarrow \infty}\Big(\int_I\int_{\Omega_{\eta(t)}}|\bu|^2\ dxdt-\int_I\int_{\Omega_{\eta_n(t)}}|\bu_n|^2\ dxdt\Big)\le 0.\]
Thus, it suffices to show that
\begin{equation}\label{eqn:kompfin}
\limsup_{n\rightarrow \infty}\Big(\int_I\int_{\Omega_{\eta_n(t)}}|\bu_n|^2\ dxdt-\int_I\int_{\Omega_{\eta(t)}}|\bu|^2\ dxdt\Big)\le 0. 
\end{equation}
While we can prove \eqref{eqn:l2konv1}$_2$ exactly as before we are not able to show \eqref{eqn:l2konv1}$_1$. This is due to the fact that the first part of \cite[Lemma A.13]{b101} is not applicable anymore. Nevertheless, defining\footnote{See Lemma \ref{lemma:mittelwertorth} for the definition of the operators $\M^\perp_{\eta_n}$.}
\begin{align*}
c_{b,n}(t)&:=\int_{\Omega_{\eta_n(t)}}\bu_n(t,\cdot)\cdot(\F_{\eta_n}\M_{\eta_n}^{\perp}b)(t,\cdot)\ dx+\int_M\pa_t\eta_n(t,\cdot)\,(\M_{\eta_n}^{\perp}b)(t,\cdot)\ dA,
\end{align*}
for $t\in\bar I$, defining $c_b$ analogously with $(\eta_n,\bu_n)$ replaced by $(\eta,\bu)$, and defining $h_n$ as in \eqref{eqn:hn}, we can make use of Lemma \ref{lemma:mittelwertorth} to show as before that $(h_n)$ tends to zero in $C(\bar I)$. An application of Lemma \ref{lemma:ehrling} yields that for
\[g_n(t):=\sup_{\norm{b}_{L^4(M)}\le
1} \big(c_{b,n}(t)-c_{b}(t)\big)\]
estimate \eqref{eqn:gest} holds for all $6/5<r<p$, thus proving that $(g_n)$ tends to zero in $L^1(I)$. Of course, we can not proceed as in the case $p>3/2$ by setting $b=\pa_t\eta_n(t,\cdot)$ since we have no bound of $\pa_t\eta_n(t,\cdot)$ in $L^4(M)$. Instead, we replace $b$ by suitable spatial-high-frequency cut-offs of the shell velocities. To this end, we fix some orthonormal basis of $L^2(M)$ and denote by $\P_k$ the orthogonal projection onto the first $k$ basis functions. By adding a zero sum, for fixed $k\in\setN$ we obtain the identity
\begin{align}\notag
&  \int_I\int_{\Omega_{\eta_n(t)}}\bu_n\cdot \F_{\eta_n}\M_{\eta_n}^{\perp}\P_k\pa_t\eta_n\ dxdt
+\int_I\im\pa_t\eta_n \M_{\eta_n}^{\perp}\mathcal{P}_k\pa_t\eta_n\ dAdt\\\notag
& -\int_I\int_{\Omega_{\eta(t)}}\bu\cdot \F_{\eta}\M_{\eta}^{\perp}\P_k\pa_t\eta\ dxdt
-\int_I\im\pa_t\eta \M_{\eta}^{\perp}\P_k\pa_t\eta\ dAdt\\\label{eqn:wursti}
&=\int_I\int_{\Omega_{\eta_n(t)}}\bu_n\cdot \F_{\eta_n}\M_{\eta_n}^{\perp}\P_k\pa_t\eta_n\ dxdt
+\int_I\im\pa_t\eta_n \M_{\eta_n}^{\perp}\P_k\pa_t\eta_n\ dAdt\\\notag
&\quad -\int_I\int_{\Omega_{\eta(t)}}\bu\cdot \F_{\eta}\M_\eta^{\perp} \P_k\pa_t\eta_n\ dxdt
-\int_I\im\pa_t\eta \M_\eta^{\perp} \P_k\pa_t\eta_n\ dAdt\\\notag
&\quad+ \int_I\int_{\Omega_{\eta(t)}}\bu\cdot
\F_{\eta}\M_{\eta}^{\perp}\P_k(\pa_t\eta_n-\pa_t\eta)\ dxdt
+ \int_I\im\pa_t\eta\,\M_\eta^{\perp} \P_k(\pa_t\eta_n-\pa_t\eta)\ dAdt.
 \end{align}
Of course, it's not a restriction to assume that the basis functions lie in $L^4(M)$.\footnote{In fact, it's this property that guarantees that the projections $\P_k$ cut-off high frequencies (in a weak, but sufficient sense).} Thus, by \eqref{ab:beschraenkt}, for fixed $k$ the first two lines of the right-hand side of \eqref{eqn:wursti} are bounded by $c\norm{g_n}_{L^1(I)}$ for some constant $c>0$. Since the sequences $(\M_\eta^{\perp} \P_k(
\pa_t\eta_n-\pa_t\eta))_n$ and $(\F_\eta\M_\eta^{\perp} \P_k(
\pa_t\eta_n-\pa_t\eta))_n$ converge to zero weakly in $L^2(I\times M)$ and $L^2(I\times B_\alpha)$, respectively, for fixed $k$ the right-hand side of \eqref{eqn:wursti} vanishes in the limit $n\rightarrow\infty$. Moreover, by adding a zero sum, we obtain
\begin{align}\notag
&  \int_I\int_{\Omega_{\eta_n(t)}}|\bu_n|^2\ dxdt+\int_I\im |\P_k\pa_t\eta_n|^2\ dAdt-\int_I\int_{\Omega_{\eta(t)}}|\bu|^2\ dxdt-\int_I\im|\P_k\pa_t\eta|^2\ dAdt\\\notag
&= \int_I\int_{\Omega_{\eta_n(t)}}\bu_n\cdot(\bu_n+\F_{\eta_n}(\M_{\eta_n}^\perp \P_k\pa_t\eta_n-\pa_t\eta_n))\ dxdt + \int_I\im\pa_t\eta_n \M_{\eta_n}^\perp \P_k\pa_t\eta_n\ dAdt\\\label{eqn:wursti2}
& -\int_I\int_{\Omega_{\eta(t)}}\bu\cdot(\bu+\F_{\eta}(\M_{\eta}^\perp \P_k\pa_t\eta-\pa_t\eta))\ dxdt - \int_I\im\pa_t\eta \M_{\eta}^\perp \P_k\pa_t\eta\ dAdt\\\notag
&-\int_I\int_{\Omega_{\eta_n(t)}}\bu_n\cdot\F_{\eta_n}\M_{\eta_n}^\perp (\P_k\pa_t\eta_n-\pa_t\eta_n)\ dxdt +\int_I\int_{\Omega_{\eta(t)}}\bu\cdot\F_{\eta}\M_{\eta}^\perp (\P_k\pa_t\eta-\pa_t\eta)\ dxdt.
 \end{align}
Here, we used the orthogonality of the projections $M_{\eta_n}^\perp$, $M_\eta^\perp$, and $\P_k$. In view of \eqref{eqn:l2konv1}$_2$ and the convergence of \eqref{eqn:wursti}, for fixed $k$ the first two lines of the right-hand side of \eqref{eqn:wursti2} vanish in the limit $n\rightarrow \infty$. Furthermore, by the definition of $\F_\eta$, see \cite{b101}, we have
\begin{align*}
\int_I&\int_{\Omega_{\eta(t)}\cap S_\alpha}\bu\cdot\F_{\eta}\M_\eta^\perp(\P_k\pa_t\eta-\pa_t\eta)\ dxdt\\
&=\int_I\int_M\int_{-\alpha}^{\eta}\exp\Big(\int_{\eta}^{s}\beta(q+\tau\,\bnu))\ d\tau\Big)\,\bnu\cdot\bu(q+s\bnu)\,|\det d\Lambda|\ ds\\
&\hspace{7cm} \M_\eta^\perp(\P_k\pa_t\eta-\pa_t\eta)\ dA(q)dt\\
&=:\int_I\int_M\psi_0\,\M_\eta^\perp(\P_k\pa_t\eta-\pa_t\eta)\ dAdt\\
&\le c\norm{\psi_0}_{L^{p}(I,W^{1,r}(M))}\norm{\M_\eta^\perp(\P_k\pa_t\eta-\pa_t\eta)}_{L^{p'}(I,(W^{1,r}(M))')}\\
&\le c\norm{\psi_0}_{L^{p}(I,W^{1,r}(M))}\norm{\M_\eta^\perp(\P_k\pa_t\eta-\pa_t\eta)}_{L^{p'}(I,(H^{1/3}(M))')}
% &\quad +\int_I\int_{\Omega\setminus\overline{S_\alpha}}\bu\cdot\F_{\eta}(\P_k\pa_t\eta-\pa_t\eta)\ dxdt
\end{align*}
for all $6/5\le r<p$. A simple calculation using Corollary \ref{lemma:spur} shows that we can bound the $L^{p}(I,W^{1,r}(M))$-norm of $\psi_0$ by the $L^{p}(I,W^{1,r}(\Omega_{\eta(t)}))$-norm of $\bu$. Moreover, we have
\begin{align*}
&\int_I\int_{\Omega\setminus\overline{S_\alpha}}\bu\cdot\F_{\eta}\M_\eta^\perp(\P_k\pa_t\eta-\pa_t\eta)\ dxdt\\
&\hspace{4cm}\le \norm{\bu}_{L^\infty(I,L^2(\Omega_{eta(t)}))}\norm{\F_{\eta}\M_\eta^\perp(\P_k\pa_t\eta-\pa_t\eta)}_{L^\infty(I,L^2(\Omega\setminus\overline{S_\alpha}))}.
\end{align*}
Remember that in $\Omega\setminus\overline{S_\alpha}$ the extension $\F_{\eta}(\P_k\pa_t\eta-\pa_t\eta)$ is given by the solution of the Stokes system with vanishing right-hand side and boundary values on $\pa(\Omega\setminus\overline{S_\alpha})$ given by
\[\exp\Big(\int_{\eta\circ q}^{-\alpha}
\beta(q+\tau\,\bnu\circ q)\ d\tau\Big)\,
(\M_\eta^\perp(\P_k\pa_t\eta-\pa_t\eta)\,\bnu)\circ q=:\psi_1\,\M_\eta^\perp(\P_k\pa_t\eta-\pa_t\eta)\circ q.\]
By a change of variables and the regularity of $\psi_1$, it's easy to see that the $(H^{1/2}(\pa(\Omega\setminus\overline{S_\alpha})))'$-norm of this function can be bounded by the $(H^{1/2}(M))'$-norm of $\M_\eta^\perp(\P_k\pa_t\eta-\pa_t\eta)$. On the other hand, Theorem 3 in \cite{b19} shows that the solution operator of the Stokes system is bounded from the space of functionals $\bg\in(H^{1/2}(\pa(\Omega\setminus\overline{S_\alpha})))'$ with $\langle \bg,\bnu\rangle=0$ to $L^2(\Omega\setminus\overline{S_\alpha})$.\footnote{Here, $\bnu$ denotes the (outer) unit normal of $\pa(\Omega\setminus\overline{S_\alpha})$.} Combining these estimates we obtain that\footnote{At first sight, it might seem awkward that we need spatial regularity of $\bu$ to control the integral over $\Omega_{\eta}\cap S_\alpha$ while this is not the case for the integral over $\Omega\setminus\overline{S_\alpha}$. Obviously, this is due to the fact that the extension operator $\F_\eta$ is not optimal in the sense that it produces no spatial regularity in $S_\alpha$.}
\[\int_I\int_{\Omega_{\eta(t)}}\bu\cdot\F_{\eta}\M_\eta^\perp(\P_k\pa_t\eta-\pa_t\eta)\ dxdt\le c\norm{\M_\eta^\perp(\P_k\pa_t\eta-\pa_t\eta)}_{L^{\infty}(I,(H^{1/3}(M))')},\]
and, similarly, we have
\[\int_I\int_{\Omega_{\eta_n(t)}}\bu_n\cdot\F_{\eta_n}\M_{\eta_n}^\perp(\P_k\pa_t\eta_n-\pa_t\eta_n)\ dxdt\le c\norm{\M_{\eta_n}^\perp(\P_k^n\pa_t\eta_n-\pa_t\eta_n)}_{L^{\infty}(I,(H^{1/3}(M))')}.\]
Using Lemma \ref{lemma:mittelwertorth}, Lemma \ref{lemma:projection}, and duality, we can make the right-hand sides small by choosing $k$ large, independently of $n$. Thus, for each $\eps>0$ we can find some fixed large $k$ such that the $\limsup$ in $n$ of the left-hand side of \eqref{eqn:wursti2} is bounded by $\eps$. This proves \eqref{eqn:kompfin} since for fixed $k$ we have
\[\limsup_{n\rightarrow \infty}\Big(\int_I\im |\P_k\pa_t\eta|^2\ dAdt-\int_I\im |\P_k\pa_t\eta_n|^2\ dAdt\Big)\le 0.\]
\qed\smallskip

\subsection{The regularized and decoupled system}
We have to regularize (and decouple) our system. As discussed in \cite{b101} it is essential to regularize the motion of the boundary. Furthermore, for technical reasons, we want to avoid to apply the proof of strong $L^2$-compactness to the Galerkin system, i.e., to the finite-dimensional approximations. For this reason, we (slightly) regularize the explicit nonlinearities in the system. Furthermore, since we want to apply monotone operator theory to the regularized system, we have to make sure that a weak solution $(\pa_t\eta,\bu)$ possesses a (formal) time-derivative in the dual of the energy class. This is achieved by perturbing the extra stress tensor $S$ into an operator $S_{\tilde\eps}$ with a $p_0$-structure for $p_0\ge 11/5$ and by adding the term $\grad_{L^2}K(\pa_t\eta)$ to the shell equation, resulting in a ``parabolization" of the whole system.\footnote{Note that the classical limit exponent $11/5$ is, in fact, not the limit exponent in our case, due to the weak regularization of the convective term announced above.} Finally, we need the weak solutions of our regularized (and decoupled) system to be unique which is most easy to prove for $p_0\ge 4$. Thus, we set $S_{\tilde\eps}(D):=S(D)+\tilde\eps |D|^{2}D$ and $p_0:=\max(p,4)$. 

%We shall use the regularization operators $\R_\eps$ (for functions defined in $\pa\Omega$ or $\bar I\times \pa\Omega$) and the modified  initial values $\bu_0^\eps$ and $\eta_1^\eps$ constructed in \cite[Subsection 3.2]{b101}. In particular, we have $\bu_0^\epsilon\in L^2(\Omega_{\R_\epsilon\eta_0})$, $\dv\bu_0^\eps=0$,  $\trnormale\bu_0^\epsilon=\eta_1^\epsilon\,\gamma(\R_\epsilon\eta_0)$, and
%\begin{equation}\label{eqn:konv7}
%  \begin{aligned}
%    \bu_0^\epsilon\rightarrow\bu_0 \qquad &\text{ in } L^2(\setR^3),\\
%    \eta_1^\epsilon\rightarrow\eta_1 \qquad&\text{ in } L^2(M).
%  \end{aligned}
%\end{equation}

We shall use the regularization operators $\R_\eps$ constructed in \cite[Subsection 3.2]{b101}. Remember that $\R_\eps\eta_0$ approximates $\eta_0$ uniformly from above. Furthermore, we note that $\trnormaln(\bu_0-\F_{\eta_0}\eta_1)=0$. Thus, extending $\bu_0-\F_{\eta_0}\eta_1\in L^2(\Omega_{\eta_0})$ by $\bfzero$ to $\setR^3$ yields a diver\-gence-free vector field in $L^2(\setR^3)$ whose support is contained in $\Omega_{\R_\eps\eta_0}$. Let $\tilde\bu_{0}^\eps$ denote a smooth divergence-free approximation of this field whose support is contained in $\Omega_{\R_\eps\eta_0}$ as well.\footnote{given, e.g., by convolution with a mollifier kernel} Moreover, let $\eta_1^\eps$ be a smooth, $C^4$ say, approximation of $\eta_1$ satisfying\footnote{given, e.g., by applying a regularization operator similar to $\R_\eps$ followed by an application of $\M_{\R_\eps\eta_0}$} \[\int_{\pa\Omega}\eta_1^\eps\,\gamma(\R_\eps\eta_0)\ dA=0,\]
and let $\bu_0^\eps:=\tilde\bu_{0}^\eps+\F_{\R_\eps\eta_0}\eta_1^\eps$. Then we have $\bu_0^\epsilon\in C^1(\overline{\Omega_{\R_\epsilon\eta_0}})$, $\dv\bu_0^\eps=0$,  $\trnormale\bu_0^\epsilon=\eta_1^\epsilon\,\gamma(\R_\epsilon\eta_0)$. From the definition of the operator $\F$ it is not hard to see that $\chi_{\Omega_{\R_\eps\eta_0}}\F_{\R_\eps\eta_0}\eta_1^\eps$ converges to $\chi_{\Omega_{\eta_0}}\F_{\eta_0}\eta_1$ in $L^2(\setR^3)$ and thus
\begin{equation}\label{eqn:konv7}
  \begin{aligned}
  \eta_1^\epsilon&\rightarrow\eta_1 &&\text{ in } L^2(M),\\
    \chi_{\Omega_{\R_\eps\eta_0}}\bu_0^\epsilon&\rightarrow\chi_{\Omega_{\eta_0}}\bu_0 &&\text{ in } L^2(\setR^3).
  \end{aligned}
\end{equation}
In the following, let $I=(0,T)$, $T>0$, be a fixed time interval and $\delta\in C(\bar I\times \pa\Omega)$ be an arbitrary, but fixed function such that $\norm{\delta}_{L^\infty(I\times\pa\Omega)}<\kappa$ and $\delta(0,\cdot)=\eta_0$. Let $\bphi$ be a vector field defined in $\Omega_{\R\delta}^I$ and $b$ a function defined in $I\times M$. 
%We extend $\bphi$ to the whole time axis by setting $\bphi(t,\cdot)=\bphi(0,\cdot)$ for $t<0$ and $\bphi(t,\cdot)=\bphi(T,\cdot)$ for $t>T$, analogously for $b$,\marginpar{Anders!} and 
For $t\in\bar I$ we define
\begin{align*}
(\R^0_\eps\bphi)(t,\cdot)&:=\T_{\R_\eps\delta(t)}\frac{1}{\eps}\int_{t-\eps}^{t}\T_{\R_\eps\delta(s)}^{-1}\bphi(s,\cdot)\ ds,\\
(\R^1_\eps b)(t,\cdot)&:=(\det(d\Psi_{\R_\eps\delta(t)}))^{-1}\frac{1}{\eps}\int_{t-\eps}^{t}\det(d\Psi_{\R_\eps\delta(s)})\,b(s,\cdot)\ ds
%\eta_0+\int_0^t\trrds\bu_k(s,\cdot)\cdot\bnu\ ds.
\end{align*}
where we extend the integrands by $0$ to the whole time axis. We have $\trred\R^0_\eps\bphi=\R^1_\eps b$ provided that $\trred\bphi=b\bnu$. Furthermore, we note that $\R^0_\eps$ preserves the divergence-free constraint. 
% let $\theta\in C^\infty([0,\infty))$ satisfy $0\le\theta\le 1$, $\theta=1$ in $[0,1]$, $\theta=0$ in $[2,\infty)$, $0\le -\theta'\le 2$, and, for $s\ge 0$, define $\theta_\eps(s)=\theta(\eps s)$ and $\Theta_\eps(s)=\int_0^s \theta_\eps(\tau)\tau\, d\tau$. 
% We set $S_\eps(D\bu)=S(D\bu)+\eps D\bu$. Thus, $S_\eps$ possesses $p_0$-structure where $p_0=\max(p,2)$.
Let us now define our decoupled and regularized problem. 
\begin{Definition}\label{def:entp}
Let $\eps,\tilde\eps>0$. A couple $(\eta,\bu)$ is called a weak solution of the decoupled and
regularized system with datum $\delta$ in the interval
  $I$ if $\eta\in Y^I\cap H^1(I,H^2_0(M))$ with $\eta(0,\cdot)=\eta_0$, $\bu\in
  X_{\R_\eps\delta,p_0}$ with $\trred \bu =\pa_t\eta\,\bnu$, and
  \begin{align}
    % &\iot \bu(t,\cdot)\cdot\bphi(t,\cdot)\ dx 
    &- \int_I\ioretd \bu\cdot\pa_t\bphi\ dxdt - \int_I\ioretd(\R^0_\eps\bu)\otimes\bu:D\bphi\, dxdt\notag\\
    &-\frac{1}{2}\int_I\im\pa_t\eta\,
    \pa_t\mathcal{R}_\eps\delta\, b\, \gamma(\mathcal{R}_\eps\delta)\ dAdt+\frac12\int_I\im(\R^1_\eps\pa_t\eta)\,\pa_t\eta\, b\, \gamma(\mathcal{R}_\eps\delta)\ dAdt
    \label{eqn:entp} \\
    & + \int_I\ioretd S_{\tilde\eps}(D\bu):D\bphi\ dxdt -\int_I\im\pa_t\eta\,
    \pa_tb\ dAdt + 2\int_I K(\eta+\eps\pa_t\eta,b)\
    dt\notag  \\
    % +\im\pa_t\eta(t,\cdot)b(t,\cdot)\ dA
    &\hspace{0.5cm}=\int_I\ioretd \ff\cdot\bphi\ dxdt + \int_I\im g\,
    b\ dAdt
    +\int_{\Omega_{\mathcal{R}_\eps\eta_0}}\bu_0^\eps\cdot\bphi(0,\cdot)\ dx +
    \im \eta_1^\eps\, b(0,\cdot)\ dA\notag 
  \end{align}
  for all test functions $(b,\bphi)\in T_{\R_\eps\delta,p}^I$.
\end{Definition}
Concerning the space of test functions, we note that $\tilde p\ge p_0$. Thus, the term involving the modified extra stress tensor is well-defined and finite. Furthermore, note that we introduced two regularization parameters, $\tilde\eps$ for the extra stress tensor and $\eps$ for the rest. The reason is that, if we let first $\eps$ tend to zero, then the explicit nonlinearity in $\pa_t\eta$ will vanish. This way, for the second limit $\tilde\eps\searrow 0$, the restriction $p>3/2$ in Proposition \ref{lemma:komp} is irrelevant.
% \begin{Definition}\label{def:entp}
% Let $1<p<\infty$. A couple $(\eta,\bu)$ is called a weak solution of the decoupled and
%   regularized problem with the data $\delta$ in the interval
%   $I$ if $\eta\in \widetilde Y^I$ with $\eta(0,\cdot)=\eta_0$, $\bu\in
%   X_{\R\delta,p_0}^I$ with $\trrd \bu =\pa_t\eta\,\bnu$, and
%   \begin{align}
%     % &\iot \bu(t,\cdot)\cdot\bphi(t,\cdot)\ dx 
%     &- \int_I\iortd \bu\cdot\pa_t\bphi\ dxdt - \int_I\iortd\bu\otimes\bu:D\bphi\, \theta(|\bu|)\, dxdt\notag\\
%     &-\frac{1}{2}\int_I\im\pa_t\eta\,
%     \pa_t\mathcal{R}\delta\, b\, \gamma(\mathcal{R}\delta)\ dAdt+\int_I\im\Theta(\pa_t\eta)\, b\, \gamma(\mathcal{R}\delta)\ dAdt
%     \label{eqn:entp} \\
%     & + \int_I\iortd S(D\bu):D\bphi\ dxdt -\int_I\im\pa_t\eta\,
%     \pa_tb\ dAdt + 2\int_I K(\eta+\pa_t\eta,b)\
%     dt\notag  \\
%     % +\im\pa_t\eta(t,\cdot)b(t,\cdot)\ dA
%     &\hspace{0.2cm}=\int_I\iortd \ff\cdot\bphi\ dxdt + \int_I\im g\,
%     b\ dAdt
%     +\int_{\Omega_{\mathcal{R}\eta_0}}\bu_0\cdot\bphi(0,\cdot)\ dx +
%     \im \eta_1\, b(0,\cdot)\ dA\notag 
%   \end{align}
%   for all test functions $(b,\bphi)\in T_{\R\delta,p_0}^I$.
% \end{Definition}

\begin{Proposition}\label{prop:entp}
Let $\eps,\tilde\eps>0$. There exists a unique weak solution $(\eta,\bu)$ of the decoupled and
regularized system with datum $\delta$ in the interval $I$ which satisfies the estimate 
\begin{equation}\label{ab:ent}
 \begin{aligned}
 &\norm{\eta}_{Y^I}^2 + \norm{\bu}_{L^\infty(I,L^2(\Omega_{\R_\eps\delta(t)}))}^2 + \norm{D\bu}_{L^p(\Omega_{\R_\eps\delta}^I)}^p\le
c_0(T,\Omega_{\R_\eps\delta}^I,\ff,g,\bu_0^\eps,\eta_0,\eta_1^\eps).
  \end{aligned}
\end{equation}
In particular, the left-hand side is bounded independently of
$\epsilon,\tilde\eps,$ and $\delta$. Furthermore, for some constant $c>0$, we have
\begin{align}\label{ab:ente}
 \eps\norm{\pa_t\eta}_{L^2(I,H^2_0(M))}^2 + \tilde\eps\norm{D\bu}_{L^{p_0}(\Omega_{\R_\eps\delta}^I)}^{p_0}\le c.
\end{align}

%\comment{M: wieso ist das wahr, da ja  $\Omega_{\R\delta}^I$ drin steht?}
\end{Proposition}

For the sake of a better readability, for the  moment, we will suppress the parameters $\epsilon,\tilde\eps$ in the notation. In particular, $\bu_0$ and $\eta_1$ denote the regularized initial values $\bu_0^\epsilon$ and $\eta_1^\epsilon$, respectively, and $S$ denotes $S_{\tilde\eps}$. For the proof of this proposition we will need the following lemma. Let 
\[E_{\eta,\bu}(t)=\frac{1}{2} \int_{\Omega_{\R\delta(t)}} |\bu(t,\cdot)|^2\ dx + \frac{1}{2} \im
|\pa_t\eta(t,\cdot)|^2\
dA + K(\eta(t,\cdot)).\]
Note that, a-priori, this function is only defined almost everywhere.

\begin{Lemma}\label{lem:approx}
 Let $(\eta,\bu)$ be a weak solution of the decoupled and
regularized system with datum $\delta$ in the interval $I$ where the field $S(D\bu)$ in \eqref{eqn:entp} may be replaced by an arbitrary field $\xi\in L^{p_0'}(\Omega_{\R\delta}^I)$. Then, $\pa_t\eta\in C(\bar I,L^2(M))$ with $\pa_t\eta(0,\cdot)=\eta_1$, $\bu\in C(\bar I,L^2(\setR^3))$ with $\bu(0,\cdot)=\bu_0$, and for all $t\in \bar I$ we have the energy identity
\begin{equation}
\begin{aligned}\label{eqn:enid}
 E_{\eta,\bu}(t)-E_{\eta,\bu}(0)=&- \int_0^t\iortd
\xi:D\bu\ dxds - 2\int_0^t K(\pa_t\eta)\ ds \\
& + \int_0^t\iortd \ff\cdot\bu\ dxds +
\int_0^t\im g\, \pa_t\eta\
dAds.
\end{aligned}
\end{equation}
\end{Lemma}
\proof
Let\footnote{Using Remark \ref{bem:tdelta2} and the fact that $\R\delta$ is smooth, we see that classical Korn's inequality holds uniformly.}
\begin{equation*}
 \begin{aligned}
V:=\{(b,\bphi)\in L^2(I,H^2_0(M)) \times L^{p_0}(I,W^{1,p_0}_{\dv}(\Omega_{\R\delta(t)}))\ |\ \trrd\bphi=b\,\bnu\}
 \end{aligned}
\end{equation*} 
and define $\bu_k:=\R^0_{1/k}\bu+\bu_{0,k}$, $\pa_t\eta_k:=\R^1_{1/k}\pa_t\eta+\eta_{1,k}$, and $\eta_k(t,\cdot):=\eta_0+\int_0^t\pa_t\eta_k(s,\cdot)\ ds$ where
\begin{align*}
\bu_{0,k}(t,\cdot)&:=(1-kt)\,\chi_{(0,1/k)}(t)\,\T_{\R\delta(t)}\T_{\R\eta_0}^{-1}\bu_0,\\
\eta_{1,k}(t,\cdot)&:=\trrdt\bu_0^k(t,\cdot)=(1-kt)\,\chi_{(0,1/k)}(t)\,\frac{\det(d\Psi_{\R\eta_0})}{\det(d\Psi_{\R\delta(t)})}\eta_1.
\end{align*}
% \begin{align*}
%  \bu_k(t,\cdot)&=\R^0_{1/k}\bu\T_{\R\delta(t)}\frac{1}{k}\int_t^{t+\frac1k}T_{\R\delta(s)}^{-1}\bu(s,\cdot)\ ds,\\
% \eta_k(t,\cdot)&=\eta_0+\int_0^t\trrds\bu_k(s,\cdot)\cdot\bnu\ ds.
% \end{align*}
Remember that $\bu_0$ and $\eta_1$ are smooth and note that $\trrd\bu_k=\pa_t\eta_k\bnu$. We claim that $\eta_k\in Y^I\cap H^2(I,H^2_0(M))$, $\bu_k\in X^I_{\R\delta,p_0}\cap W^{1,p_0}(I,L^{p_0}(\Omega_{\R\delta(t)}))$, that
\begin{equation*}%\label{eqn:konvp1}
  (\pa_t\eta_k,\bu_k)\rightarrow(\pa_t\eta,\bu)\quad\text{ in } V
\end{equation*}
for $k\rightarrow\infty$, and that the functionals
\begin{align}\notag
\big\langle\frac{d}{dt}(\pa_t\eta_k,\bu_k),(b,\bphi)\big\rangle:= &\int_I\iortd \pa_t\bu_k\cdot\bphi\ dxdt +\frac{1}{2}\int_I\im\pa_t\eta_k\,\pa_t\mathcal{R}\delta\, b\, \gamma(\mathcal{R}\delta)\ dAdt\\\label{eqn:timefunc}
&+\int_I\im\pa_t^2\eta_k\,b\ dAdt
\end{align}
are bounded in $V'$. Except for the inclusion $(\pa_t\bu_k)\subset L^{p_0}(\Omega^I_{\R\delta})$ and the boundedness of the functionals, these claims are obvious if we remember that $\T_{\R\delta}$ is an isomorphism between the involved function spaces on $\Omega_{\R\delta}^I$ and the corresponding function spaces on $I\times\Omega$, see Remark \ref{bem:tdelta2}. Before proving the remaining assertions, let us draw the relevant conclusions. We can proceed as in \cite[Remark 1.17]{b101} to show that the extension of $\bu_k$ by $(\pa_t\eta_k\bnu)\circ q$ lies in $C(\bar I,L^{p_0}(B_\alpha))$ for $\norm{\R\delta}_{L^\infty(I\times M)}<\alpha<\kappa$. Thus, the extension of $\bu_k$ by $\boldsymbol{0}$ lies in $C(\bar I,L^2(\setR^3))$. For all $s,t\in [0,T]$ we have
\begin{equation}\label{eqn:zeitabl}
 \begin{aligned}
\big\langle\frac{d}{dt}(\pa_t\eta_k,\bu_k),(\pa_t\eta_k,\bu_k)\chi_
{(s,t)}\big\rangle=&\ \frac{1}{2} \int_{\Omega_{\R\delta(t)}} |\bu_k(t,\cdot)|^2\ dx +
\frac{1}{2} \im |\pa_t\eta_k(t,\cdot)|^2\
dA\\
&-\frac{1}{2} \int_{\Omega_{\R\delta(s)}} |\bu_k(s,\cdot)|^2\ dx -
\frac{1}{2} \im |\pa_t\eta_k(s,\cdot)|^2\ dA.
 \end{aligned}
\end{equation}
Replacing $\eta_k$ by $\eta_k-\eta_l$ and $\bu_k$ by $\bu_k-\bu_l$ and integrating the resulting identity over $I$ with respect to $s$, we obtain
\begin{multline*}
\int_{\Omega_{\R\delta(t)}} |(\bu_k-\bu_l)(t,\cdot)|^2\ dx + \im
|\pa_t(\eta_k-\eta_l)(t,\cdot)|^2\ dA \\
\shoveright{\le c\,\big(\bignorm{\frac{d}{dt}(\pa_t(\eta_k-\eta_l),(\bu_k-\bu_l))}_{V'}\norm{(\pa_t(\eta_k-\eta_l),
\bu_k-\bu_l)
}_V +\norm{\bu_k-\bu_l}_{L^2(\Omega_{\R\delta}^I)}^2}\\
+\norm{\pa_t\eta_k-\pa_t\eta_l}_{L^2(I\times M)}^2\big).
 \end{multline*}
Extending the functions $\bu_k$ and $\bu$ by $\boldsymbol{0}$ to $I\times \setR^3$, we deduce from this estimate and the properties of the approximations that the sequences $(\bu_k)$
and $(\pa_t\eta_k)$ converge to $\bu$ in $C(\bar I, L^2(\setR^3))$ and to $(\pa_t\eta)$ in $C(\bar I,
L^2(M))$, respectively. By an argument analogous to the one given in \cite[Remark 3.3]{b101}, using the $L^2$-continuity of $\pa_t\eta$ and $\bu$, we can show that
\eqref{eqn:entp} holds with $\bu(0,\cdot)$ und $\pa_t\eta(0,\cdot)$ in place of $\bu_0$ and $\eta_1$, respectively, proving that $\pa_t\eta(0,\cdot)=\eta_1$, $\bu(0,\cdot)=\bu_0$.

Choosing $s=0$ in \eqref{eqn:zeitabl}, the right-hand side converges to
\begin{equation*}
 \begin{aligned}
 \frac{1}{2} \int_{\Omega_{\R\delta(t)}} |\bu(t,\cdot)|^2\ dx +
\frac{1}{2} \im |\pa_t\eta(t,\cdot)|^2\
dA\\
-\frac{1}{2} \int_{\Omega_{\R\delta(0)}} |\bu_0|^2\ dx -
\frac{1}{2} \im |\eta_1|^2\
dA.
 \end{aligned}
\end{equation*}
By the uniform boundedness of $(d/dt(\pa_t\eta_k,\bu_k))$ in $V'$, the left-hand side converges to $\langle\Sigma,(\pa_t\eta,\bu)\chi_{(0,t)}\rangle$ where $\Sigma\in V'$ is given by
 \begin{align}
  \langle\Sigma,(b,\bphi)\rangle_V=&  \int_I\iortd(\R^0\bu)\otimes\bu:D\bphi\, dxdt-\frac12\int_I\im(\R^1\pa_t\eta)\,\pa_t\eta\, b\, \gamma(\mathcal{R}\delta)\ dAdt\notag\\
    & - \int_I\iortd \chi:D\bphi\ dxdt - 2\int_I K(\eta+\pa_t\eta,b)\ dt\notag +\int_I\iortd \ff\cdot\bphi\ dxdt\\
    & + \int_I\im g\, b\ dAdt\label{eqn:sigma}
 \end{align}
for $(b,\bphi)\in V$. This can be seen as follows. For $(b,\bphi)\in T_{\R\delta,p}^I$ with $b(0,\cdot)=0$, $\bphi(0,\cdot)=\boldsymbol{0}$, an application of Reynold's transport theorem shows that
\begin{align*}
\big\langle\frac{d}{dt}(\pa_t\eta_k,\bu_k),(b,\bphi)\big\rangle= &-\int_I\iortd \bu_k\cdot\pa_t\bphi\ dxdt -\frac{1}{2}\int_I\im\pa_t\eta_k\,\pa_t\mathcal{R}\delta\, b\, \gamma(\mathcal{R}\delta)\ dAdt\\
&-\int_I\im\pa_t\eta_k\,\pa_t b\ dAdt.
% -\int_{\Omega_{\mathcal{R}\eta_0}}\bu_k(0,\cdot)\cdot\bphi(0,\cdot)\ dx\\
% & - \im \pa_t\eta_k(0,\cdot)\, b(0,\cdot)\ dA.
\end{align*}
Now, if we let $k\rightarrow\infty$, use \eqref{eqn:entp}, and note that these test functions are dense in $V$,\footnote{In order to prove the denseness, we can proceed analogously to the proof of denseness employed in \cite{b101} just before the proof of Proposition 3.15.} we obtain \eqref{eqn:sigma}. Noting that
\[(\R^0\bu)\otimes\bu:D\bu=(\R^0\bu)\cdot\nabla\frac{|\bu|}{2},\]
we see that for $(b,\bphi)=(\pa_t\eta,\bu)\chi_{(0,t)}$ the first two terms on the right-hand side of \eqref{eqn:sigma} cancel. Hence, \eqref{eqn:enid} follows.

Now, let us prove the inclusion $(\pa_t\bu_k)\subset L^{p_0}(\Omega^I_{\R\delta})$ and the boundedness of the functionals \eqref{eqn:timefunc}. Note that
\begin{equation}
\begin{aligned}\label{eqn:tderiv}
 \pa_t(\R^0_{1/k}\bu_k)(t,\cdot)=&\bigg(d\bigg(\frac{d\Psi_{\R\delta(t)}}{\det d\Psi_{\R\delta(t)}}k\int_{t-\frac1k}^{t}T_{\R\delta(s)}^{-1}\bu(s,\cdot)\ ds\bigg)\bigg)\circ \Psi_{\R\delta(t)}^{-1}\ \pa_t\Psi_{\R\delta(t)}^{-1}\\
&+\bigg(\pa_t\bigg(\frac{d\Psi_{\R\delta(t)}}{\det d\Psi_{\R\delta(t)}}k\int_{t-\frac1k}^{t}T_{\R\delta(s)}^{-1}\bu(s,\cdot)\ ds\bigg)\bigg)\circ \Psi_{\R\delta(t)}^{-1}.
\end{aligned}
\end{equation}
While the second term on the right-hand side obviously lies in $L^{p_0}(\Omega^I_{\R\delta})$, the same is true for the first term because the spatial derivatives of $\bu$ lie in $L^{p_0}(\Omega^I_{\R\delta})$. We can similarly show that $\pa_t\bu_{0,k}\in L^{p_0}(\Omega^I_{\R\delta})$. Let us proceed with the boundedness of the functionals \eqref{eqn:timefunc}. For $(b,\bphi)\in T_{\R\delta,p}^I$ we have
\begin{align}\notag
\big\langle\frac{d}{dt}(\pa_t\eta_k,\bu_k),&(b,\bphi)\big\rangle= \int_I\int_{\Omega_{\R\delta(t)}}\pa_t\bu_{0,k}\cdot\bphi\ dxdt
+\int_I\int_M\pa_t^2\eta_{1,k}\,b\ dAdt\\\label{eqn:horstsepp}
&-\int_I\iortd \R_{1/k}^0\bu\cdot\pa_t\bphi\ dxdt -\frac{1}{2}\int_I\im \pa_t\eta_k\,\pa_t\mathcal{R}\delta\, b\, \gamma(\mathcal{R}\delta)\ dAdt\\\notag
&-\int_I\im(\R^1_{1/k}\pa_t\eta)\,\pa_t b\ dAdt.
%-\int_{\Omega_{\R\eta_0}}\bu_0\cdot\bphi(0,\cdot)\ dx-\int_M\eta_1\,b(0,\cdot)\ dA.
% -\int_{\Omega_{\mathcal{R}\eta_0}}\bu_k(0,\cdot)\cdot\bphi(0,\cdot)\ dx\\
% & - \im \pa_t\eta_k(0,\cdot)\, b(0,\cdot)\ dA.
\end{align}
In order to deal with the third and the fifth term on the right hand side, let us denote the $L^2(\Omega^I_{\R\delta})$-adjoints of $\R^0_{1/k}$ and $\R^1_{1/k}$ by $(\R^0_{1/k})'$ and $(\R^1_{1/k})'$, respectively. We have $(\R^1_{1/k})'=\R^1_{-1/k}$ and
\[((\R^0_{1/k})'\bphi)(t,\cdot)=(\T_{\R\delta(t)}^{-1})'k\int_{t}^{t+\frac1k}(\T_{\R\delta(s)})'\bphi(s,\cdot)\ ds\]
where we extend the integrands by $0$ to the whole time axis and
\begin{align*}
 (\T_{\R\delta(t)})'\bphi(t,\cdot)=(d\Psi_{\R\delta(t)})^T\bphi(t,\cdot)\circ\Psi_{\R\delta(t)},\\
(\T_{\R\delta(t)}^{-1})'\widetilde\bphi(t,\cdot)=(d\Psi_{\R\delta(t)}^{-1})^T\widetilde\bphi(t,\cdot)\circ\Psi_{\R\delta(t)}^{-1}.
\end{align*}
We compute
\begin{align*}
 -\int_I\iortd \R_{1/k}^0&\bu\cdot\pa_t\bphi\ dxdt=-\int_I\iortd \bu\cdot(\R^0_{1/k})'\pa_t\bphi\ dxdt\\
&=-\int_I\iortd \bu\cdot\pa_t((\R^0_{1/k})'\bphi)\ dxdt+\int_I\iortd \bu\cdot [\pa_t,(\R^0_{1/k})']\bphi\ dxdt.
\end{align*}
Here, the commutator
\begin{equation}
\begin{aligned}\label{eqn:tder}
([\pa_t,(\R^0_{1/k})']\bphi)(t,\cdot)&=[\pa_t,(\T_{\R\delta(t)}^{-1})']k\int_{t}^{t+\frac1k}(\T_{\R\delta(s)})'\bphi(s,\cdot)\ ds\\
&\quad + (\T_{\R\delta(t)}^{-1})'k\int_{t}^{t+\frac1k}[\pa_s,(\T_{\R\delta(s)})']\bphi(s,\cdot)\ ds.
\end{aligned}
\end{equation}
is acting derivatively only on the spatial variable of $\bphi$, cf. \eqref{eqn:tderiv}. Thus, the $L^{p_0}(\Omega^I_{\R\delta})$-norm of $[\pa_t,(\R^0_{1/k})']\bphi$ is bounded by the $L^{p_0}(I,W^{1,p_0}(\Omega_{\R\delta(t)}))$-norm of $\bphi$. Analogously, we have
\begin{align*}
 -\int_I\im (\R_{1/k}^1\pa_t\eta)\,\pa_t b&\ dAdt=-\int_I\im \pa_t\eta\, (\R^1_{1/k})'\pa_t b\ dAdt\\
&=-\int_I\im \pa_t\eta \,\pa_t((\R^1_{1/k})'b)\ dAdt +\int_I\im \pa_t\eta\, [\pa_t,\R^1_{1/k})']b\ dAdt,
\end{align*}
and, here, the $L^2(I\times M)$-norm of $[\pa_t,\R^1_{1/k})']b$ is even bounded by the $L^2(I\times M)$-norm of $b$. Unfortunately, in general, $((\R^0_{1/k})'\bphi,(\R^1_{1/k})'b)\notin T^I_{\R\delta,p}$ since the adjoint operator $(\R^0_{1/k})'$ preserves neither the divergence-free constraint nor the structure of the boundary values. We can overcome this problem by replacing $((\R^0_{1/k})'\bphi,(\R^1_{1/k})'b)$ by $(\R^0_{-1/k}\bphi,\R^1_{-1/k} b)\in T^I_{\R\delta,p}$. In order to do so, remembering that $(\R^1_{1/k})'=\R^1_{-1/k}$, we need to bound the functionals
\begin{align}\label{eqn:func}
(b,\bphi)\mapsto &\int_I\iortd \bu\cdot\pa_t(((\R^0_{1/k})'-\R^0_{-1/k})\bphi)\ dxdt
%&-\int_I\im \pa_t\eta \,\pa_t(((\R^1_{1/k})'-\R^1_{-1/k}) b)\ dAdt
\end{align}
in $V'$. As we saw above the commutators $[\pa_t,\T_{\R\delta}]$ and $[\pa_t,(\T_{\R\delta}^{-1})']$ are acting derivatively only on the spatial variable so that the corresponding terms in \eqref{eqn:func} can be estimated by a constant multiple of
%L^{p_0}(\Omega^I_{\R\delta})$-norm of $[\pa_t,((\R^0_{1/k})'-\R^0_{-1/k})]\bphi$ can be bounded by
the $L^{p_0}(I,W^{1,p_0}(\Omega_{\R\delta(t)}))$-norm of $\bphi$. Hence, we need to deal with the terms resulting from the time-derivatve acting on the Steklov means. These terms evaluated at $t\in I$ give
\begin{align*}
&k\,(\T_{\R\delta(t)}^{-1})'\Big((\T_{\R\delta(t+\frac1k)})'\bphi(t+1/k,\cdot)-(\T_{\R\delta(t)})'\bphi(t,\cdot)\Big)\notag\\
& -k\,\T_{\R\delta(t)}\Big(\T_{\R\delta(t+\frac1k)}^{-1}\bphi(t+1/k,\cdot)-\T_{\R\delta(t)}^{-1}\bphi(t,\cdot)\Big)\\
&\quad=k\,\Big(\T_{\R\delta(t)}\T_{\R\delta(t+\frac1k)}^{-1}-(\T_{\R\delta(t)}^{-1})'(\T_{\R\delta(t+\frac1k)})'\Big)\bphi(t+1/k,\cdot).
\end{align*}
But from the smoothness of $\R\delta$ we can deduce that
\begin{align}
&\int_I\int_{\Omega_{\R\delta(t)}}\big|\Big(\T_{\R\delta(t)}\T_{\R\delta(t+\frac1k)}^{-1}-(\T_{\R\delta(t)}^{-1})'(\T_{\R\delta(t+\frac1k)})'\Big)\bphi(t+1/k,\cdot)\big|^{p_0}\ dxdt\notag\\\label{eq:horst2}
&=\int_I\int_\Omega \big|d\Psi_{\R\delta(t)}(\det d\Psi_{\R\delta(t)})^{-1}\big(d\Psi_{\R\delta(t+\frac1k)}^{-1}(\det d\Psi_{\R\delta(t+\frac1k)}^{-1})^{-1}\big)\circ\Psi_{\R\delta(t+\frac1k)}\\
&\quad-(d\Psi_{\R\delta(t)}^{-1})^T\circ\Psi_{\R\delta(t)}\, d\Psi_{\R\delta(t+\frac1k)}^T\big|^{p_0}\ |\bphi(t+1/k,\cdot)\circ\Psi_{\R\delta(t+\frac1k)}|^{p_0}\ \det d\Psi_{\R\delta(t)} dxdt\notag\\\notag
&\le \frac{c}{k}\int_I\int_\Omega |\bphi(t+1/k,\cdot)\circ\Psi_{\R\delta(t+\frac1k)}|^{p_0}\ dxdt\le \frac{c}{k}
\end{align}
proving the boundedness of \eqref{eqn:func} in $V'$. In view of \eqref{eqn:horstsepp} and from what we showed so far it remains to bound the functionals
\begin{align}\notag
(b,\bphi)&\mapsto \int_I\int_{\Omega_{\R\delta(t)}}\pa_t\bu_{0,k}\cdot\bphi\ dxdt
+\int_I\int_M\pa_t^2\eta_{1,k}\,b\ dAdt-\int_I\iortd \bu\cdot\pa_t(\R_{-1/k}^0\bphi)\ dxdt\\\label{eqn:funci}
& -\frac{1}{2}\int_I\im \pa_t\eta\,\pa_t\mathcal{R}\delta\, (\R^1_{-1/k} b)\, \gamma(\mathcal{R}\delta)\ dAdt
-\int_I\im\pa_t\eta\,\pa_t(\R^1_{-1/k} b)\ dAdt.
% -\int_{\Omega_{\mathcal{R}\eta_0}}\bu_k(0,\cdot)\cdot\bphi(0,\cdot)\ dx\\
% & - \im \pa_t\eta_k(0,\cdot)\, b(0,\cdot)\ dA.
\end{align}
in $V'$. Using \eqref{eqn:entp} we can replace the last three terms in \eqref{eqn:funci} by 
\[\int_{\Omega_{\mathcal{R}\eta_0}}\bu_0\cdot(\R^0_{-1/k}\bphi)(0,\cdot)\ dx +
    \im \eta_1\, (\R^1_{-1/k}b)(0,\cdot)\ dA\]
since the remaining terms in \eqref{eqn:entp} can be bounded by a constant multiple of $\norm{(b,\bphi)}_{V}$, reflecting the fact that the formal time-derivative of $(\pa_t\eta,\bu)$ lies in $V'$. Here, the convective term is the crucial one. By interpolation, we have
\begin{align*}
 \int_I\iortd(\R^0\bu)&\otimes\bu:D(\R_{-1/k}^0\bphi)\, dxdt\\
&\le \norm{\R^0\bu}_{L^{11/3}(\Omega_{\R\delta}^I)}\norm{\bu}_{L^{11/3}(\Omega_{\R\delta}^I)}\norm{D\R_{-1/k}^0\bphi}_{L^{11/5}(\Omega_{\R\delta}^I)}\\
&\le c\norm{\bu}_{L^{\infty}(I,L^2(\Omega_{\R\delta(t)})}^{4/5}\norm{\bu}_{L^{11/5}(I,W^{1,11/5}(\Omega_{\R\delta(t)})}^{6/5}\norm{\bphi}_{L^{11/5}(I,W^{1,11/5}(\Omega_{\R\delta(t)})}.
\end{align*}
Thus, the lemma is proved if we can bound the functionals
\begin{equation}
\begin{aligned}
(b,\bphi)\mapsto &\int_I\int_{\Omega_{\R\delta(t)}}\pa_t\bu_{0,k}\cdot\bphi\ dxdt
+\int_I\int_M\pa_t^2\eta_{1,k}\,b\ dAdt\\\label{eqn:seppl}
&+\int_{\Omega_{\mathcal{R}\eta_0}}\bu_0\cdot(\R^0_{-1/k}\bphi)(0,\cdot)\ dx +
    \im \eta_1\, (\R^1_{-1/k}b)(0,\cdot)\ dA
% -\int_{\Omega_{\mathcal{R}\eta_0}}\bu_k(0,\cdot)\cdot\bphi(0,\cdot)\ dx\\
% & - \im \pa_t\eta_k(0,\cdot)\, b(0,\cdot)\ dA.
\end{aligned}
\end{equation}
in $V'$. A simple computation shows that the sum of the first and the third term equals
\begin{align*}
k\int_0^{\frac1k}\int_{\Omega_{\R\delta(t)}}\big((\T_{\R\delta(t)}^{-1})'\T_{\R\eta_0}'-\T_{\R\delta(t)}\T_{\R\eta_0}^{-1}\big)\bu_0\cdot\bphi\ dxdt\\
+\int_0^{\frac1k}\int_{\Omega_{\R\delta(t)}}(1-kt)\,\pa_t\big(\T_{\R\delta(t)}\T_{\R\eta_0}^{-1}\bu_0\big)\cdot\bphi\ dxdt.
\end{align*}
This expression can be bounded by a constant multiple of the $L^{p_0}(\Omega_{\R\delta}^I)$-norm of $\bphi$ analogously to \eqref{eq:horst2} and \eqref{eqn:tderiv}. The sum of the second and the fourth term in \eqref{eqn:seppl} can be handled similarly. This completes the proof.
\qed\medskip

\proof (of Proposition \ref{prop:entp}) We use
the Galerkin method. We proceed exactly as in \cite{b101} for the construction of time-dependent basis functions $(\bW_k)$ and $W_k$ such that 
\begin{align*}
  \spann\{(\phi\, W_k,\phi\,\bW_k)\ |\ \phi\in C_0^1([0,T)),\,
  k\in\setN\}
\end{align*}
is dense in $T^I_{\R\delta,p_0}$. We seek
functions $\alpha_n^k:[0,T]\rightarrow\setR$, $k,n\in\setN$, such that
$\bu_n:=\alpha_n^k\,\bW_k$ and $\eta_n(t,\cdot):=\int_0^t\alpha^k_n\,
W_k\ ds+\eta_0$ (summation with respect to $k$ from $1$ to $n$) solve
the equation%\footnote{As usual we suppress
%  the dependence of the functions on $t$.}
\begin{align}\notag
    % &\iot \bu(t,\cdot)\cdot\bphi(t,\cdot)\ dx
    &\iortd \pa_t\bu_n\cdot\bW_j\ dx - \iortd (\R^0\bu_n)\otimes\bu_n:D\bW_j\ dx+\frac{1}{2}\int_I\im\pa_t\eta\,
    \pa_t\mathcal{R}\delta\, b\, \gamma(\mathcal{R}\delta)\ dA   \\\notag &\quad+\frac12\im\R^1(\pa_t\eta_n)\,\pa_t\eta_n\, W_j\, \gamma(\mathcal{R}\delta)\ dA
 +\iortd S(D\bu_n):D\bW_j\ dx+\im\pa^2_t\eta_n\, W_j\ dA\\ 
 &\hspace{3cm}+ 2K(\eta_n+\pa_t\eta_n, W_j)\label{eqn:gal}
    % +\im\pa_t\eta(t,\cdot)b(t,\cdot)\ dA 
    =\iortd \ff_n\cdot\bW_j\ dx + \im g_n\, W_j\ dA
 \end{align}
for all $1\le j \le n$. Here, $\ff_n$ and $g_n$ are smooth functions
which converge to $\ff$ and $g$ in $L^2_\loc([0,\infty)\times
\setR^3)$ and $L^2_\loc([0,\infty)\times M)$, respectively. As in \cite{b101}, we construct initial conditions $\alpha^k_n(0)$ such that
\begin{align*}
  \pa_t\eta_n(0,\cdot)&\rightarrow\eta_1\qquad \text{ in } L^2(M),
  \\
  \bu_n(0,\cdot)&\rightarrow \bu_0\qquad \text{ in }L^2(\Omega_{\R\eta_0}).
\end{align*}
With these initial conditions, \eqref{eqn:gal} is a
Cauchy problem for a system of ordinary
integro-differential equations of the form ($1\le j\le n$, summation
with respect to $k$ from $1$ to $n$)
\begin{equation*}
  A_{jk}(t)\,\dot\alpha^k(t) = B_{j}(t,\alpha(t)) + \int_0^t
  C_{j}(\alpha(t),\alpha(s),t,s)\ ds +   D_j(t).
\end{equation*}
Here, the functions $A_{jk},D_j:[0,T]\rightarrow\setR$ and $B_j:[0,T]\times\setR^n\rightarrow\setR$, given by
\begin{align*}
    A_{jk}(t)&=\iortd \bW_k\cdot\bW_j\ dx + \im W_k\, W_j\ dA,\\
    B_{j}(t,\alpha(t))&=\Big(-\iortd \pa_t\bW_k\cdot\bW_j\ dx - \im \pa_t W_k\, W_j\ dA-2K(W_k,W_j)\\
&\quad\ \ \ - \frac12\im W_k\,W_j\,\gamma(\R\delta)\ dA\Big)\alpha^k(t) - \iortd S(\alpha^k(t)\,D\bW_k):D\bW_j\ dx,\\
    D_j(t)&=\iortd \ff_n\cdot \bW_j\ dx + \im g_n\, W_j\ dA,
 \end{align*}
are continuous, while the functions $C_j:\setR^d\times\setR^d\times [0,T]\times [0,T]\rightarrow\setR$, given by
\[C_{j}(\alpha(t),\alpha(s),t,s)=-2K(W_k(s),W_j(t))\,\alpha^k(s)+\frac{1}{\eps}\chi_{(t-\eps,t)}(s)\,E_{jkl}(t)\,\alpha^k(s)\,\alpha^l(t)\]
with
\[E_{jkl}(t)=\ \iortd\bW_l\otimes\bW_k:D\bW_j\ dx -\frac12\im W_l\,W_k\,W_j\,\gamma(\R\delta)\ dA,\]
are measurable and bounded on compact subsets of their domain. Furthermore, we saw in \cite{b101} that the matrices $A(t)$ are invertible. Now, one can easily adapt the proof of Peano's existence theorem to show that there exists a unique, local $C^1$-solution $\alpha$ which exists as long as $|\alpha(t)|$ stays bounded, cf. \cite[Appendix A.3]{b62}. Let us now test \eqref{eqn:gal} with $(\pa_t\eta_n,\bu_n)$. We saw in the proof of Lemma \ref{lem:approx} that the second and the fourth term on the left-hand side cancel, while the first and the third term yield
\[\frac{d}{dt}\frac12\int_{\Omega_{\R\delta(t)}}|\bu_n|^2\ dx.\]
Thus, we can procceed as in Subsection 1.3 to obtain
\begin{align*}
  \norm{\eta_n}_{Y^I}^2 + \norm{\bu_n}_{L^\infty(I,L^2(\Omega_{\R\delta(t)}))}^2 + \norm{D\bu_n}_{L^p(\Omega_{\R\delta}^I)}^{p}\le  c_0(T,\Omega_{\R\delta}^I,\ff_n,g_n,\bu_n(0,\cdot),\eta_0,\pa_t\eta_n(0,\cdot)),
\end{align*}
as well as
\begin{align*}
\eps\norm{\pa_t\eta_n}_{L^2(I,H^2_0(M))} + \tilde\eps\norm{D\bu_n}_{L^{p_0}(\Omega_{\R_\eps\delta}^I)}\le c
\end{align*} 
for some constant $c>0$. In particular, the solutions exist on the whole time interval $[0,T]$. From these bounds and \eqref{eqn:tderiv} we deduce that
\begin{align}\label{ab:reg}
\norm{S(D\bu_n)}_{L^{p_0'}(\Omega_{\R\delta}^I)}+\norm{\R^0\bu_n}_{W^{1,p_0}(\Omega_{\R\delta}^I)}+\norm{\R^1\pa_t\eta_n}_{H^1(I\times M)}\le c'
\end{align}
for another constant $c'>0$. Hence, for a subsequence (again denoted by the index $n$) we have
  \begin{alignat*}{2}
    \eta_n&\rightarrow \eta &&\quad\text{weakly in } H^1(I,H^2_0(M)),\\ 
\pa_t\eta_n&\rightarrow \pa_t\eta &&\quad\text{weakly$^*$ in } L^\infty(I,L^2(M)),\\ 
\R^1\pa_t\eta_n&\rightarrow \R^1\pa_t\eta &&\quad\text{ in } L^2(I\times M),\\ 
    \bu_n&\rightarrow\bu &&\quad\text{ weakly$^*$ in }
    X_{\R\delta,p_0}^I,\\
\R^0\bu_n&\rightarrow\R^0\bu &&\quad\text{ in } L^{p_0}(\Omega_{\R\delta}^I),\\
S(D\bu_n)&\rightarrow\xi &&\quad\text{ weakly in }
    L^{p_0'}(\Omega_{\R\delta}^I).
  \end{alignat*}
The above convergences and $\trrd \bu_n=\pa_t\eta_n\,\bnu$ imply the identity $\trrd
\bu=\pa_t\eta\,\bnu$. Furthermore, by the lower semi-continuity of the norms with respect to weak and weak* convergence, we deduce \eqref{ab:ent} and \eqref{ab:ente}. Multiplying \eqref{eqn:gal} by $\phi(t)$, where $\phi\in C_0^1([0,T))$, integrating over $I$, integrating by parts in time, letting $n\rightarrow\infty$, and using the denseness of the test functions in $T_{\R\delta,p}^I$,\footnote{The denseness can be shown by exactly the same argument used in \cite{b101} just before the proof of Proposition 3.15.} we see that the couple $(\eta,\bu)$ satisfies \eqref{eqn:entp} with $\xi$ in place of $S(D\bu)$. Thus, it remains to identify $\xi$. In view of the $L^2$-continuity of $\pa_t\eta$ and $\bu$, we can proceed analogously to the proof of \cite[Remark 3.3]{b101} to show that, for $(b,\bphi)\in T_{\R\delta,p}^I$ with the constraint $b(T,\cdot)=0$, $\bphi(T,\cdot)=\boldsymbol{0}$ replaced by $b(0,\cdot)=0$, $\bphi(0,\cdot)=\boldsymbol{0}$, $(\eta,\bu)$ satifies \eqref{eqn:entp} with $\xi$ in place of $S(D\bu)$ and with the right-hand side replaced by
\[\int_I\iortd \ff\cdot\bphi\ dxdt + \int_I\im g\,
    b\ dAdt
    -\int_{\Omega_{\mathcal{R}\eta_0}}\bu(T,\cdot)\cdot\bphi(T,\cdot)\ dx -
    \im \eta(T,\cdot)\, b(T,\cdot)\ dA.\]
On the other hand, we note that subsequences of $(\bu_n(T,\cdot))$ and $(\pa_t\eta_n(T,\cdot))$ converge weakly to functions $\bu^*$ and $\eta^*$ in $L^2(\Omega_{\R\delta(T)})$ and $L^2(M)$, respectively. Thus, multiplying $\eqref{eqn:gal}$ by $\phi(t)$, $\phi\in C_0^1((0,T])$, and taking the limit as before, we see that, for $(b,\bphi)$ as above, $(\eta,\bu)$ satifies \eqref{eqn:entp} with $\xi$ in place of $S(D\bu)$ and with the right-hand side replaced by
\[\int_I\iortd \ff\cdot\bphi\ dxdt + \int_I\im g\,
    b\ dAdt
    -\int_{\Omega_{\mathcal{R}\eta_0}}\bu^*\cdot\bphi(T,\cdot)\ dx -
    \im \eta^*\, b(T,\cdot)\ dA.\]
This yields $\bu^*=\bu(T,\cdot)$, $\eta^*=\pa_t\eta(T,\cdot)$. Furthermore, a subsequence of  $(\eta_n(T,\cdot))$ converges to $\eta(T,\cdot)$ weakly in $H^2_0(M)$. We have already seen that the Galerkin solutions satisfy the energy identity
\begin{equation*}
 \begin{aligned}
&\int_I\int_{\Omega_{\R\delta(t)}}S(D\bu_n):D\bu_n\ dxdt + 2\int_I K(\pa_t\eta_n)\ dt\\ 
&\hspace{1cm}=-E_{\eta_n,\bu_n}(T)+E_{\eta_n,\bu_n}(0) + \int_I\iortd \ff_n\cdot\bu_n\ dxdt
+\int_I\im g_n\,\pa_t\eta_n\ dAdt.
 \end{aligned}
\end{equation*}
Taking the $\limsup$ of this equation, eploiting the weak lower semi-continuity of the energy $E$,\footnote{Note that each continuous, non-negative quadratic form, e.g. $K$, is weakly lower semi-continuous. This follows by taking the $\liminf$ of the inequality
\[0\le K(\eta_n-\eta,\eta_n-\eta)=K(\eta_n)-2K(\eta,\eta_n)+K(\eta).\]} and noting that $\eta_n(0,\cdot)=\eta_0$ for all $n\in\setN$, we obtain
\begin{equation*}
 \begin{aligned}
&\limsup_n \int_I\int_{\Omega_{\R\delta(t)}}S(D\bu_n):D\bu_n\ dxdt + 2\int_I K(\pa_t\eta_n)\ dt\\ 
&\hspace{1cm}\le-E_{\eta,\bu}(T)+E_{\eta,\bu}(0) + \int_I\iortd \ff\cdot\bu\ dxdt +\int_I\im g\,
\pa_t\eta\ dAdt.
 \end{aligned}
\end{equation*}
From the energy identity \eqref{eqn:enid} for the weak solution $(\eta,\bu)$ (with $\xi$ in place of $S(D\bu)$) we deduce that
\begin{equation*}
 \begin{aligned}
&\limsup_n \int_I\int_{\Omega_{\R\delta(t)}}S(D\bu_n):D\bu_n\ dxdt + 2\int_I K(\pa_t\eta_n)\ dt\\ 
&\hspace{4cm}\le\int_I\int_{\Omega_{\R\delta(t)}}\xi:D\bu\ dxdt + 2\int_I K(\pa_t\eta)\ dt.
 \end{aligned}
\end{equation*} 
Using this estimate and the weak convergences, we obtain
\begin{equation*}
 \begin{aligned}
  0&\le\limsup_n\Big(\int_I\int_{\Omega_{\R\delta(t)}}(S(D\bu_n)-S(D\bu)):(D\bu_n-D\bu)\ dxdt\\
&\hspace{4.3cm}+2\int_I K(\pa_t\eta_n-\pa_t\eta,\pa_t\eta_n-\pa_t\eta)\ dt\Big)\\
&=\limsup_n\Big(\int_I\int_{\Omega_{\R\delta(t)}}S(D\bu_n):D\bu_n + S(D\bu):D\bu\\
&\hspace{4.75cm} - S(D\bu_n):D\bu -
S(D\bu):D\bu_n\ dxdt\\
&\hspace{2cm}+2\int_I K(\pa_t\eta_n,\pa_t\eta_n) + K(\pa_t\eta,\pa_t\eta)-
2K(\pa_t\eta_n,\pa_t\eta)\ dt\Big)\\
&\le \int_I\int_{\Omega_{\R\delta(t)}}\xi:D\bu + S(D\bu):D\bu - \xi:D\bu - S(D\bu):D\bu\ dxdt\\
&\hspace{2.2cm}+2\int_I K(\pa_t\eta,\pa_t\eta) + K(\pa_t\eta,\pa_t\eta)-
2K(\pa_t\eta,\pa_t\eta)\ dt\\
&=0.
 \end{aligned}
\end{equation*}
Hence, for a subsequence, we have
\[(S(D\bu_n)-S(D\bu)):(D\bu_n-D\bu)\rightarrow 0\]
a.e. in $\Omega_{\R\delta}^I$. By Proposition \ref{lemma:dmm} we infer that
$D\bu_n\rightarrow D\bu$ and hence $S(D\bu_n)\rightarrow S(D\bu)$ a.e. in
$\Omega_{\R\delta}^I$. Finaly, Vitali's convergence theorem yields
$\xi=S(D\bu)$. This proves the existence of weak solutions. 

Now, let us show uniqueness. For weak solutions $(\eta_0,\bu_0)$ and $(\eta_1,\bu_1)$ of the regularized and decoupled system with datum $\delta$ in the interval $I$ their difference $(\eta:=\eta_1-\eta_0,\bu:=\bu_1-\bu_0)$ is a weak solution, too, with $S(D\bu)$ replaced by $\xi=S(D\bu_1)-S(D\bu_0)$ and $\ff$, $g$ replaced by
\begin{align*}
\tilde\ff:=-(\R^0\bu\cdot\nabla)\bu_0 - (\R^0\bu_0\cdot\nabla)\bu,\quad\tilde g:=\frac12 (\R^1\pa_t\eta)\,\pa_t\eta_0 + \frac12 (\R^1\pa_t\eta_0)\,\pa_t\eta,
\end{align*}
respectively. Thus, by Lemma \ref{lem:approx}, we have
\begin{align*}
E_{\eta,\bu}(t)&\le -  \int_0^t\iortd (\R^0\bu\cdot\nabla)\bu_0\cdot\bu\ dxds +
\int_0^t\im \frac12 (\R^1\pa_t\eta)\,\pa_t\eta_0\, \pa_t\eta\ dAds\\
&\le \int_0^t \norm{\R^0\bu(s)}_{L^2(\Omega_{\R\delta(s)})}\norm{\bu(s)}_{L^2(\Omega_{\R\delta(s)})}\norm{\nabla\bu_0(s)}_{L^4(\Omega_{\R\delta(s)})}\ ds\\
&\quad+\int_0^t\norm{\R^1\pa_t\eta(s)}_{L^2(M)}\norm{\pa_t\eta(s)}_{L^2(M)}\norm{\pa_t\eta_0(s)}_{L^4(M)}\ ds\\
&\le \int_0^t\big(\norm{\bu(s)}_{L^2(\Omega_{\R\delta(s)})}^2+\norm{\pa_t\eta(s)}_{L^2(M)}^2\big)\,c(s)\ ds
\end{align*}
for some nonnegative function $c\in L^1(I)$. In the first inequality we used the monotonicity of $S$ and the fact that the second terms in $\tilde\ff$ and $\tilde g$ cancel when tested against $(\pa_t\eta,\bu)$. By Gronwall's inequality, we have $E_{\eta,\bu}\equiv 0$, proving that the solutions coincide.
\qed

\subsection{Fixed-point argument}
Let us now define solutions of our regularized problem. 
\begin{Definition} \label{def:epsilon} 
Let $\eps,\tilde\eps>0$. A couple $(\eta,\bu)$ is a weak solution of the $(\eps,\tilde\eps)$-regularized system in the interval $I$ if $\eta\in Y^I\cap H^1(I,H^2_0(M))$ with $\norm{\eta}_{L^\infty(I\times M)}<\kappa$, $\eta(0,\cdot)=\eta_0$,
  and $\bu\in X_{\R_\eps\eta,p_0}^I$ with $\trre\bu = \pa_t\eta\,\bnu$, and
  \begin{align}
    % &\iot \bu(t,\cdot)\cdot\bphi(t,\cdot)\ dx 
    &- \int_I\iorte \bu\cdot\pa_t\bphi\ dxdt - \int_I\iorte(\R_\eps^0\bu)\otimes\bu:D\bphi\, dxdt\notag\\
    &-\frac{1}{2}\int_I\im\pa_t\eta\,
    \pa_t\R_\eps\eta\, b\, \gamma(\R_\eps\eta)\ dAdt+\frac12\int_I\im\R_\eps^1(\pa_t\eta)\,\pa_t\eta\, b\, \gamma(\R_\eps\delta)\ dAdt
    \label{eqn:entp0} \\
    & + \int_I\iorte S_{\tilde\eps}(D\bu):D\bphi\ dxdt -\int_I\im\pa_t\eta\,
    \pa_t b\ dAdt + 2\int_I K(\eta+\eps\pa_t\eta,b)\
    dt\notag  \\
    % +\im\pa_t\eta(t,\cdot)b(t,\cdot)\ dA
    &\hspace{0.2cm}=\int_I\iorte \ff\cdot\bphi\ dxdt + \int_I\im g\,
    b\ dAdt
    +\int_{\Omega_{\R_\eps\eta_0}}\bu_0\cdot\bphi(0,\cdot)\ dx +
    \im \eta_1\, b(0,\cdot)\ dA\notag 
  \end{align}
  for all test functions $(b,\bphi)\in T_{\R_\eps\eta,p}^I$.
\end{Definition}
% \begin{Remark}
%   Note that due to the definition \eqref{eqn:taueta} we get that
%   $\tau(\eta_\epsilon)< \infty$. If we additionally require that
%   $\norm{\eta_\epsilon}_{L^\infty(I\times M)}\le \alpha$ for some
%   $\alpha<\kappa$ then we get ${\tau(\eta_\epsilon)\le c(\alpha)}$
%   uniformly with respect to $\epsilon$.
% \end{Remark}
% \vspace{0.2cm}
%\smallskip

\begin{Proposition}\label{lemma:approx}
  There exists a $T>0$ such that for all sufficiently small $\epsilon,\tilde\eps>0$
  there exists a weak solution $(\eta,\bu)$ of the $(\eps,\tilde\eps)$-regularized system in the interval $I=(0,T)$. Furthermore, we have
  \begin{equation}\label{ab:approx}
    \begin{aligned}
      &\norm{\eta}_{Y^I}^2 + \norm{\bu}_{L^\infty(I,L^2(\Omega_{\R_\eps\eta(t)}))}^2 + \norm{D\bu}_{L^p(\Omega_{\R_\eps\eta}^I)}^p\le
      c_0(T,\Omega_{\R_\epsilon\eta}^I,\ff,g,\bu^\epsilon_0,\eta_0, 
      \eta^\epsilon_1) 
    \end{aligned}
  \end{equation}
  and $\sup_{\epsilon}\tau(\eta_\epsilon)<\infty$. The time $T$ can be
  chosen to depend only on $\tau(\eta_0)$ and the bound \eqref{ab:approx} for the $Y^I$-norm of $\eta_\epsilon$. Finally, for some constant $c>0$, we have
\begin{align}\label{ab:approxe}
 \eps\norm{\pa_t\eta}_{L^2(I,H^2_0(M))}^2 + \tilde\eps\norm{D\bu}_{L^{p_0}(\Omega_{\R_\eps\eta}^I)}^{p_0}\le c.
\end{align}
\end{Proposition}
%\comment{M: wir sollten zum beweis von
%  $\sup_{\epsilon}\tau(\eta_\epsilon)<\infty$ etwas sagen.}
\proof
We set $\alpha:=(\norm{\eta_0}_{L^\infty(M)}+\kappa)/2$ and
fix arbitrary but sufficiently small  $\epsilon,\tilde\eps>0$. For a better
readability, in the following, we will omit the symbols $\epsilon,\tilde\eps$. We want to use Schauder's fixed point theorem. To this end, we define the space $Z:=C(\bar I\times\pa\Omega)$ with the closed, convex subset 
\begin{align*}
  D:=\big\{(\delta,\bv)\in Z\ |\ \delta(0,\cdot)=\eta_0,\,
  \norm{\delta}_{L^\infty(I\times\pa\Omega)}\le\alpha\big\}.
\end{align*}
Let $F:D\rightarrow Z$ map each $\delta\in D$ to the component $\eta$ of
the unique weak solution $(\eta,\bu)$ of the decoupled and regularized
system with datum $\delta$. From \eqref{ab:ent} we deduce that the norm of $\eta$ in
\begin{align}\label{eq:emb}
  Y^I\embedding C^{0,1-\theta}(\bar I, C^{0,2\theta
    -1}(\pa\Omega))\quad (1/2<\theta<1)
\end{align}
is bounded. Since $\eta(0,\cdot)=\eta_0$, we can choose the time
interval $I=(0,T)$ so small that $\norm{\eta}_{L^\infty(I\times M)}\le\alpha$, independently of the parameter $\epsilon$; in particular, $\tau(\eta)\le c(\alpha)$. Thus,
$F$ maps $D$ into itself. Furthermore, from \eqref{eq:emb} and the compact embedding of the H\"older space into $Z$, we see that $F(D)$ is relatively compact in $Z$. It remains to show that $F$ is continuous. To this end, we let $(\delta_n)\subset D$ be a sequence converging to $\delta$ in $Z$ and $(\eta_n,\bu_n)$ be the corresponding weak solutions given by Proposition \ref{lem:approx}. 
% The proof of the relative compactness of $(\pa_t\eta_n,\bu_n)$ in $L^2(I\times M)\times L^2(I\times\setR^3)$ can be taken almost literally from the proof of Proposition \eqref{lemma:komp}. The only differences are the slightly changed form \eqref{eqn:entp} of the system and that in some places one has to replace the sequence $(\eta_n)$ with limit $\eta$ by the sequence $(\R\delta_n)$ with limit $\R\delta$. Moreover, due to the regularization one can, in fact, simplify some of the
% arguments. 
In view of \eqref{ab:ent}, \eqref{ab:ente}, and \eqref{ab:reg}, we can deduce that for a subsequence we have
  \begin{alignat}{2}\label{eqn:wichtig}\notag
    \eta_n&\rightarrow \eta &&\quad\text{ weakly in } H^1(I,H^2_0(M))\text{ and, thus, uniformly},\\\notag 
\pa_t\eta_n&\rightarrow \pa_t\eta &&\quad\text{
      weakly$^*$ in } L^\infty(I,L^2(M)),\\ \notag
\R^1\pa_t\eta_n&\rightarrow \R^1\pa_t\eta &&\quad\text{ in } L^2(I\times M),\\
    \bu_n&\rightarrow\bu &&\quad\text{ weakly$^*$ in }
    L^\infty(I,L^2(\setR^3)),\\\notag
\nabla\bu_n&\rightarrow\nabla\bu &&\quad\text{ weakly$^*$ in }
    L^{p_0}(I\times\setR^3)),\\\notag
\R^0\bu_n&\rightarrow\R^0\bu &&\quad\text{ in } L^{p_0}(I\times\setR^3),\\\notag
S(D\bu_n)&\rightarrow\xi &&\quad\text{ weakly in }
    L^{p_0'}(\Omega_{\R\delta}^I).
  \end{alignat}
Here, we extend $\nabla\bu_n$, $S(D\bu_n)$, and $\nabla\bu$, which are a-priori defined on
$\Omega_{\R\delta_n}^I$ and $\Omega_{\R\delta}^I$, respectively, by $\bfzero$
to the whole of $I\times \setR^3$. We have to show that $\eta=F(\delta)$. The property $\eta(0,\cdot)=\eta_0$ follows immediately from the uniform convergence of $(\eta_n)$. Moreover, we can show exactly like in the proof of \cite[Proposition 3.35]{b101} that
$\pa_t\eta\,\bnu=\trrd\bu$. It remains to prove that
\eqref{eqn:entp} is satisfied.  For all $n$ and all test
functions $(b_n,\bphi_n)\in T_{\R\delta_n,p}^I$, we have
 \begin{align}
    % &\iot \bu(t,\cdot)\cdot\bphi(t,\cdot)\ dx 
    &- \int_I\iortdn \bu_n\cdot\pa_t\bphi_n\ dxdt - \int_I\iortdn(\R^0\bu_n)\otimes\bu_n:D\bphi_n\, dxdt\notag\\
    &-\frac{1}{2}\int_I\im\pa_t\eta_n\,
    \pa_t\mathcal{R}\delta_n\, b_n\, \gamma(\mathcal{R}\delta_n)\ dAdt+\frac12\int_I\im\R^1(\pa_t\eta_n)\,\pa_t\eta_n\, b_n\, \gamma(\mathcal{R}\delta_n)\ dAdt
    \label{eqn:entp2} \\
    & + \int_I\iortdn S(D\bu_n):D\bphi_n\ dxdt -\int_I\im\pa_t\eta_n\,
    \pa_tb_n\ dAdt + 2\int_I K(\eta_n+\pa_t\eta_n,b_n)\
    dt\notag  \\
    % +\im\pa_t\eta(t,\cdot)b(t,\cdot)\ dA
    &\hspace{0.5cm}=\int_I\iortdn \ff\cdot\bphi\ dxdt + \int_I\im g\,
    b\ dAdt
    +\int_{\Omega_{\mathcal{R}\eta_0}}\bu_0\cdot\bphi_n(0,\cdot)\ dx +
    \im \eta_1\, b_n(0,\cdot)\ dA\notag.
  \end{align}
As in \cite{b101}, we can pass to the limit in this equation (with the exception of the extra stress tensor) by using, for given $(b,\bphi)\in T_{\R\delta_n,p}^I$, the special test functions $(b_n,\bphi_n):=(\M_{\R\delta_n}b,\F_{\R\delta_n}\M_{\R\delta_n}b)\in
  T_{\R\delta_n,p}^I$  on the one hand, and the test functions $(0,\bphi)\in T_{\R\delta}^I$ with $\bphi(T,\cdot)=0$ and
  $\supp\bphi\subset\Omega_{\R\delta}^{\bar I}$ on the other hand. Here, additionally to \eqref{eqn:wichtig}, we need to take into account the assertions $\rm (1.b)$, $\rm (2.b)$ in Lemma \ref{lemma:konvergenzen}. Finally, in order to identify the function $\xi$, we can proceed almost literally as in the proof of Proposition \eqref{lem:approx}. Essentially, we only have to replace the integrals over $\Omega_{\R\delta}^I$ by integrals over $I\times \setR^3$, extending the corresponding functions by $0$ to $I\times \setR^3$. This shows that $(\eta,\bu)$ is the unique weak solution of the decoupled and regularized
problem with datum $\delta$. Thus, $F$ is continuous and, by Schauder's fixed point theorem, possesses a fixed point. This concludes the proof.\qed

\smallskip

\subsection{Limiting process}
Now, we can prove our main result by letting the
regularizing parameters in Definition \ref{def:epsilon} tend
to zero.\\

\beweis (of Theorem \ref{theorem:hs}) First, we fix $\tilde\eps>0$ and let $\eps\searrow 0$. We have shown that there exists
a $T>0$ such that for all $\epsilon=1/n$, $n\in\setN$ sufficiently
large, there exists a weak solution $(\eta_\epsilon,\bu_\epsilon)$ of
the $(\eps,\tilde\eps)$-regularized problem in the interval
$I=(0,T)$. For a subsequence, the estimates \eqref{ab:approx}, \eqref{ab:approxe}, Proposition \ref{prop:korn}, and the compact embedding $Y^I\compactembedding C(\bar I\times\pa\Omega)$ yield the following
convergences
  \begin{alignat}{2}\notag
    \eta_\epsilon,\,\R_\epsilon\eta_\epsilon&\rightarrow\eta &&\quad\text{
      weakly$^*$ in }  L^\infty(I,H^2_0(M))\text{ and uniformly},\\\notag
    \pa_t\eta_\epsilon,\,
    \pa_t\R_\epsilon\eta_\epsilon&\rightarrow\pa_t\eta &&\quad\text{
      weakly$^*$ in }
    L^\infty(I,L^2(M)),\\\label{eqn:konv3}
    \bu_\epsilon&\rightarrow\bu &&\quad\text{ weakly$^*$ in }
    L^\infty(I,L^2(\setR^3)),\\\notag
S(D\bu_\eps)&\rightarrow\xi &&\quad\text{ weakly in }
    L^{p_0'}(I\times\setR^3).
\end{alignat}
Here, we extend $D\bu_\eps$, $\nabla\bu_\eps$, $S(D\bu_\eps)$, and $D\bu$ which are a-priori defined only on
$\Omega_{\R_\eps\eta}^I$ and $\Omega_{\eta}^I$, respectively, by $0$ to the whole of $I\times \setR^3$. The uniform convergence of $(\R_\epsilon\eta_\epsilon)$ follows from the estimate 
\begin{align*}
  \norm{ \R_\epsilon\eta_\epsilon-\eta}_{L^\infty(I\times \pa\Omega)}
  \le\norm{\R_\epsilon(\eta_\epsilon-\eta)}_{L^\infty(I\times
    \pa\Omega)}+\norm{\R_\epsilon\eta-\eta}_{L^\infty(I\times
    \pa\Omega)}.
\end{align*}
Now, we can repeat the proof of  Proposition \ref{lemma:komp} almost literally to show that
\begin{equation}\label{eqn:konv31}
  \begin{aligned}
    \pa_t\eta_\epsilon&\rightarrow\pa_t\eta &&\text{ in }
    L^2(I\times M),\\
    \bu_\epsilon,\R^0_\eps\bu_\eps&\rightarrow\bu &&\text{ in }
    L^2(I\times\setR^3).
  \end{aligned}
\end{equation}
From these convergences and the definition of $\R^0_\eps$, $\R^1_\eps$ it is not hard to see that
\begin{equation}\label{eqn:konv32}
  \begin{aligned}
    \R^1_\eps\pa_t\eta_\epsilon&\rightarrow\pa_t\eta &&\text{ in }
    L^2(I\times M),\\
    \R^0_\eps\bu_\eps&\rightarrow\bu &&\text{ in }
    L^2(I\times\setR^3).
  \end{aligned}
\end{equation}
As in the proof \cite[Proposition 3.35]{b101}, we obtain the
identity $\tr\bu=\pa_t\eta\,\bnu$. By \eqref{ab:approx}, Proposition \ref{prop:korn}, and Corollary \ref{lemma:sobolev}, the sequence $(\bu_n)$ is bounded in $L^\infty(I,L^2(\setR^3))\cap L^{p}(I,L^r(\setR^3))$ for all $r<p$.\footnote{In fact, this is true for $p_0$ in place of $p$. But since we will have to repeat this argument when taking the limit $\tilde\eps\searrow 0$, we don't want to make use of this fact.} From this bound, \eqref{eqn:konv31}, and \eqref{eqn:konv32}, by interpolation, we obtain that
\begin{equation}\label{eqn:konv20}
  \begin{aligned}
    &\bu_\epsilon,\R^0_\eps\bu_\eps\rightarrow\bu&&\text{ in }L^r(I\times \setR^3)
  \end{aligned}
\end{equation}
for all $1\le r<10p/6$. Similarly, we deduce from \eqref{ab:approx}, \eqref{eqn:konv31}, and \eqref{eqn:konv32} that
\begin{equation}\label{eqn:konv21}
  \begin{aligned}
    &\pa_t\eta_\epsilon,\R^1_\eps\pa_t\eta_\eps,\pa_t\R_\eps\eta_\eps\rightarrow\pa_t\eta&&\text{ in }L^4(I,L^2(M)).
  \end{aligned}
\end{equation}
  The lower semi-continuity of the norms yields the estimate
\eqref{ab:hs}, while the uniform convergence of $(\eta_\epsilon)$
gives $\eta(0,\cdot)=\eta_0$. For all $\epsilon=1/n$, $n$ sufficiently large, and all
$(b_\epsilon,\bphi_\epsilon)\in T_{\R_\epsilon\eta_\epsilon,p}^I$ we
have
\begin{align}\label{eqn:final}
    % &\iot \bu(t,\cdot)\cdot\bphi(t,\cdot)\ dx 
    &- \int_I\iorte \bu_\eps\cdot\pa_t\bphi_\eps\ dxdt - \int_I\iorte(\R_\eps^0\bu_\eps)\otimes\bu_\eps:D\bphi_\eps\, dxdt\\\notag
    &-\frac{1}{2}\int_I\im\pa_t\eta_\eps\,
    \pa_t\R_\eps\eta_\eps\, b_\eps\, \gamma(\R_\eps\eta_\eps)\ dAdt+\frac12\int_I\im\R_\eps^1(\pa_t\eta_\eps)\,\pa_t\eta_\eps\, b_\eps\, \gamma(\R_\eps\eta_\eps)\ dAdt
    \label{eqn:entpeps} \\
    & + \int_I\iorte S_{\tilde\eps}(D\bu_\eps):D\bphi_\eps\ dxdt -\int_I\im\pa_t\eta_\eps\,
    \pa_t b_\eps\ dAdt + 2\int_I K(\eta_\eps+\eps\pa_t\eta_\eps,b_\eps)\
    dt\notag  \\
    % +\im\pa_t\eta(t,\cdot)b(t,\cdot)\ dA
    &\hspace{0.2cm}=\int_I\iorte \ff\cdot\bphi_\eps\ dxdt + \int_I\im g\,
    b_\eps\ dAdt
    +\int_{\Omega_{\R_\eps\eta_0}}\bu_0^\eps\cdot\bphi_\eps(0,\cdot)\ dx +
    \im \eta_1^\eps\, b_\eps(0,\cdot)\ dA.\notag 
  \end{align}
Note that, by \eqref{ab:approxe}, we have
\begin{align*}
 \big|\int_I K(\eps\pa_t\eta_\eps,b_\eps)\ dt\big|\le \eps\,c\norm{\pa_t\eta_\eps}_{L^2(I,H^2_0(M))}\norm{b_\eps}_{L^2(I,H^2_0(M))}\le \sqrt{\eps}\,c\norm{b_\eps}_{L^2(I,H^2_0(M))},
\end{align*}
and thus, this term vanishes in the limit.  Just like in \cite{b101}, we make use of the special test functions
  $(b_\epsilon,\bphi_\epsilon):=(\M_{\R_\epsilon\eta_\epsilon}b,
  \F_{\R_\epsilon\eta_\epsilon}\M_{ \R_\epsilon\eta_\epsilon}b)\in
  T_{\R_\epsilon\eta_\epsilon,p}^I$ for given $(b,\bphi)\in T_{\eta,p}^I$ on the one hand, and the test functions $(0,\bphi)\in T_{\eta,p}^I$ with $\bphi(T,\cdot)=0$
  and $\supp\bphi\subset\Omega^{\bar I}_{\eta}$ on the other hand. Using \eqref{eqn:konv3}, \eqref{eqn:konv20},
  \eqref{eqn:konv21}, and \eqref{eqn:konv7}, as well as the assertions $\rm
  (1.b)$, $\rm (2.b)$ in Lemma \ref{lemma:konvergenzen}, we can pass to the limit in \eqref{eqn:final} with the exception of the extra stress tensor. The
  convergences \eqref{eqn:konv20} are needed for the second term, while the convergences \eqref{eqn:konv21} are needed for the third and the fourth term. This implies the validity of \eqref{eqn:schwach} with $\eps=0$ and $S(D\bu)$ replaced by $\xi$ for $(b,\bphi)=(b,\F_{\eta}b)\in T_{\eta,p}^I$; note that the third and the fourth term in \eqref{eqn:entp0} cancel for $\eps=0$.

Now, let us identify $\xi$. Let us fix some open, bounded interval $J$ and some open ball $B$ such that for the cylinder $Q:=J\times B$ we have $\overline Q\subset \Omega_\eta^I$. By the uniform convergence of $\R_\eps\eta_\eps$, we have $\overline Q\subset \Omega_{\R_\eps\eta_\eps}^I$ for sufficiently small $\eps$. Thus, setting $\bv_\eps:=\bu_\eps-\bu$ and $G_\eps:=G^0_\eps+G^1_\eps$ with
\begin{align*}
 G^0_\eps&:=S_{\tilde\eps}(D\bu_\eps)-\xi,\\
G^1_\eps&:=\R^0_\eps\bu_\eps\otimes\bu_\eps-\bu\otimes\bu,
\end{align*}
for a subsequence we have\footnote{Note that classical Korn's inequality holds in $B$.}
\begin{alignat*}{2}
 \bv_\eps&\rightarrow \boldsymbol{0} &&\quad\text{ weakly in }
    L^{p_0}(J,W^{1,p_0}(Q)),\\
G^0_\eps&\rightarrow 0 &&\quad\text{ weakly in }
    L^{p_0'}(Q),\\
G^1_\eps&\rightarrow 0 &&\quad\text{ in }
    L^{r}(Q)\text{ for all }1\le r<5p/6
\end{alignat*}
and
\begin{align*}
 -\int_Q\bv_\eps\cdot\pa_t\bphi+G_\eps:\nabla\bphi\ dxdt=0
\end{align*}
for all vector fields $\bphi\in C_0^\infty(I\times B)$ with $\dv\bphi=0$ and sufficiently small $\eps$. By H\"older's inequality and Theorem \ref{thm:lipschitz} {\it (a)}, we deduce that for $\theta\in (0,1)$ and $\zeta\in C_0^\infty(\frac16 Q)$ with $\chi_{\frac18 Q}\le\zeta\le\chi_{\frac16 Q}$ we have
\begin{align*}
&\limsup_{\eps\searrow 0}\int_Q\big((S_{\tilde\eps}(D\bu_\eps)-S_{\tilde\eps}(D\bu)):D(\bu_\eps-\bu)\big)^\theta\zeta\ dxdt\\
&\qquad\le c2^{(\theta-1)k} + \limsup_{\eps\searrow 0}\int_Q (S_{\tilde\eps}(D\bu_\eps)-S_{\tilde\eps}(D\bu)):\nabla\bv_\eps\,\zeta\,\chi_{O_{n,k}^c}\ dxdt.
\end{align*}
The second term on the right-hand side is bounded by
\begin{align*}
 \limsup_{\eps\searrow 0}\int_Q G^0_\eps:\nabla\bv_\eps\,\zeta\,\chi_{O_{n,k}^c}\ dxdt+\limsup_{\eps\searrow 0}\int_Q(\xi-S_{\tilde\eps}(D\bu)):\nabla\bv_\eps\,\zeta\,\chi_{O_{n,k}^c}\ dxdt.
\end{align*}
Here, by Theorem \ref{thm:lipschitz} {\it (b)}, the first term is bounded by $c2^{-k/p}$, while, by Theorem \ref{thm:lipschitz} {\it (a)} and the weak convergence of $(\nabla\bv_\eps)$, the second term can be estimated by
\begin{align*}
\limsup_{\eps\searrow 0}\int_Q(\xi-S_{\tilde\eps}(D\bu)):\nabla\bv_\eps\,\zeta\ dxdt
 + \limsup_{\eps\searrow 0} \int_Q(\xi-S_{\tilde\eps}(D\bu)):\nabla\bv_\eps\,\zeta\,\chi_{O_{n,k}}\ dxdt\le c2^{-k}.
\end{align*}
Thus, we showed that
\begin{align*}
&\limsup_{\eps\searrow 0}\int_Q\big((S_{\tilde\eps}(D\bu_\eps)-S_{\tilde\eps}(D\bu)):D(\bu_\eps-\bu)\big)^\theta\zeta\ dxdt=0
\end{align*}
and, hence, for a subsequence
\[(S_{\tilde\eps}(D\bu_\eps)-S_{\tilde\eps}(D\bu)):D(\bu_\eps-\bu)\rightarrow 0\]
a.e. in $\frac18 Q$. As in the proof of Proposition \ref{prop:entp} we can infer that $\xi=S_{\tilde\eps}(D\bu)$ in $\frac18 Q$. Since $Q$ was arbitrary, we have $\xi=S_{\tilde\eps}(D\bu)$ in $\Omega_\eta^I$.

We can repeat these arguments almost literally for the limit $\tilde\eps\searrow 0$. The main difference is that, in general, we won't have strong $L^2$-compactness of $(\pa_t\eta_{\tilde\eps})$ if $p\le 3/2$. But since equation \eqref{eqn:schwach} (with $S$ replaced by $S_{\tilde\eps}$) contains no (explicit) nonlinearities in $(\pa_t\eta_{\tilde\eps})$, strong compactness is not needed. Furthermore, in this limit process we have to use that
\[\int_I\int_{\Omega_{\eta_{\tilde\eps}(t)}} \tilde\eps|D\bu_{\tilde\eps}|^{p_0-2}D\bu_{\tilde\eps}:D\bphi_{\tilde\eps}\ dxdt\le \tilde\eps\norm{D\bu_{\tilde\eps}}_{L^p(\Omega_{\eta_{\tilde\eps}}^I)}^{p-1}\norm{D\bphi_{\tilde\eps}}_{L^p(\Omega_{\eta_{\tilde\eps}}^I)}\le c\tilde\eps^{1/p}.\]
Finally, we can proceed as in \cite{b101} to show that the solution exists as long as the inequality $\norm{\eta(t,\cdot)}_{L^\infty(M)}<\kappa$ holds.
  \qed

%%%%%%%%%%%%%%%%%%%%%%%%%%%%%%%%%%%%%%%%%%%%%%%%%%%%%%%%%%%%%%%%%%%%%%%%%%%%%%%%%%%%%%%%%%%%%%%%%%%%%%%%%%%%%%%%%%%%%%%%%%%%%%%%%%
%Anhang
\appendix
\section{Appendix}\label{sec:app} \label{sec:ode}

The following classical result can be found in \cite{b13}.
\begin{Proposition}{\bf (Reynolds
    transport theorem)}\label{theorem:reynolds}
  Let $\Omega\subset\setR^3$ be a bounded domain with $C^1$-boundary,
  let $I\subset\setR$ be an interval, and let $\Psi\in
  C^1(I\times\overline\Omega,\setR^3)$ such that
  \[
    \Psi_t:=\Psi(t,\cdot):\overline\Omega\rightarrow
    \Psi_t(\overline\Omega)
  \] 
  is a diffeomorphism for all $t\in I$. We set
  $\Omega_t:=\Psi_t(\Omega)$ and $\bv:=(\pa_t\Psi)\circ\Psi_t^{-1}$.
  Then we have  for  all $\xi\in C^1(\bigcup_{t\in I}\{t\}\times
  \overline\Omega_t)$ and $t\in I$ 
  \[
    \frac{d}{dt}\int_{\Omega_t} \xi(t,x)\ dx= \int_{\Omega_t} \pa_t
    \xi(t,x)\ dx + \int_{\pa\Omega_t} \bv\cdot \bnu_t\ \xi(t,\cdot)\
    dA_t.
  \]
  Here, $dA_t$ denotes the surface measure and $\bnu_t$ denotes the
  outer unit normal of $\pa\Omega_t$.
\end{Proposition}
\smallskip

\begin{Theorem}\label{thm:lipschitz}
 Let $1<p,r<\infty$, $d\in\setN$, $B\subset\setR^d$ an open ball, and $J$ an open, bounded interval. We set $Q:=J\times B$ and assume that the vector fields $\bv_n$ and the $\setR^{d\times d}$-valued fields $G^0_n$ and $G^1_n$, $n\in\setN$, satisfy
\begin{alignat*}{2}
 \bv_n&\rightarrow\boldsymbol{0}&&\quad\text{ weakly in } L^p(J,W^{1,p}(B)),\\
G^0_n&\rightarrow 0&&\quad\text{ weakly in } L^{p'}(Q),\\
G^1_n&\rightarrow 0&&\quad\text{ strongly in } L^{r}(Q).
\end{alignat*}
Furthermore, we assume that the sequence $(\bv_n)$ is bounded in $L^\infty(I,L^r(B))$ and that for $G_n=G_n^0+G_n^1$ and all vector fields $\bphi\in C_0^\infty(I\times B)$ with $\dv\bphi=0$ we have
\[\int_Q\bv_n\cdot\pa_t\bphi+G_n:\nabla\bphi\ dxdt=0.\]
Then,\footnote{For $\alpha>0$ we denote by $\alpha Q$ the cylinder $Q$ scaled by $\alpha$ with respect to its center.} for $\zeta\in C_0^\infty(\frac16 Q)$ with $\chi_{\frac18 Q}\le\zeta\le\chi_{\frac16 Q}$ and all $n,k\in\setN$, $k$ sufficiently large, there exist open sets $O_{n,k}\subset Q$ such that
\begin{itemize}
%  \item[\it (b)] $\bu_{n,k}\in L^s(J,W^{1,s}_{0,\dv}(B))$ for all $1<s<\infty$ and $\supp(u_{n,k})\subset \frac16 Q$,
% \item[\it (c)] $\bu_{n,k}=\bu_n$ a.e. in $\frac18 Q\setminus O_{n,k}$,
% \item[\it (d)] $\norm{\nabla\bu_{n,k}}_{L^\infty(I,\text{BMO}(B))}\le c\lambda_{n,k}$,
% \item[\it (e)] $\bu_{n,k}\rightarrow \boldsymbol{0}$ in $L^\infty(Q)$ for fixed $k$ and $n\rightarrow\infty$,
% \item[\it (f)] $\nabla\bu_{n,k}\rightarrow 0$ weakly in $L^s(Q)$ for all $1<s<\infty$, fixed $k$, and $n\rightarrow\infty$, 
\item[\it (a)] $\limsup_{n\rightarrow\infty}|O_{n,k}|\le c 2^{-k}$,
\item[\it (b)] $\limsup_{n\rightarrow\infty}|\int_Q G_n^0:\nabla\bu_{n}\,\zeta\,\chi_{O_{n,k}^c}\ dxdt|\le c 2^{-k/p}$.
% \item[\it (h)] $\limsup_{n\rightarrow\infty}\norm{\nabla\bu_{n,k}\chi_{O_{n,k}}}_{L^p(B)}^p\le c\lambda_{n,k}^p|O_{n,k}|$,
%\item[\it (i)] $\limsup_{n\rightarrow\infty}|\int_Q G_n:\nabla\bu_{n,k}\ dxdt|\le c\lambda_{n,k}|O_{n,k}|$.
\end{itemize}
\end{Theorem}
\proof See \cite{b109}.\qed\medskip

\begin{Proposition}\label{lemma:dmm}
Let $S:M_{sym}\rightarrow M_{sym}$ be continuous and strictly monotone, i.e., \[(S(A)-S(B)):(A-B)>0\] for $A,B\in M_{sym}$, $A\not= B$. Furthermore, let $(A_n)_{n\in\setN}\subset M_{sym}$ be a sequence such that
\[\lim_{n\rightarrow\infty}\ (S(A_n)-S(A)):(A_n-A)=0\]
for some $A\in M_{sym}$. Then $\lim_{n\rightarrow\infty}A_n=A$.
\end{Proposition}
\proof See \cite{b111}.\smallskip

\begin{Lemma}\label{lemma:mittelwert}
Let $\eta\in Y^I$ with $\norm{\eta}_{L^\infty(I\times M)}<\kappa$. There exists a linear operator
$\M_\eta$ such that
\begin{equation*}
 \begin{aligned}
  \norm{\M_\eta b}_{L^r(I\times M)}&\le c\, \norm{b}_{L^r(I\times M)},\\
  \norm{\M_\eta b}_{C(\bar I,L^r(M)}&\le c\, \norm{b}_{C(\bar I,L^r(M)},\\
  \norm{\M_\eta b}_{L^r(I,H^2_0(M))}&\le c\, \norm{b}_{L^r(I,H^2_0(M))},\\
  \norm{\M_\eta b}_{H^1(I,L^2(M))}&\le c\, \norm{b}_{H^1(I,L^2(M))}\\
 \end{aligned}
\end{equation*}
for all $1\le r\le\infty$ and
\[\int_M (\M_\eta b)(t,\cdot)\,\gamma(\eta(t,\cdot))\ dA = 0\]
for almost all $t\in I$. The constant $c$ depends only on $\Omega$, $\norm{\eta}_{Y^I}$ and
$\tau(\eta)$; it stays bounded as long as $\norm{\eta}_{Y^I}$ and $\tau(\eta)$ stay bounded.
\end{Lemma}
\proof We only give the definition of the operator, for a proof of the claims see \cite{b101}. For fixed, but arbitrary $\psi\in C_0^\infty(\inn M)$ with $\psi\ge
0$, $\psi\not\equiv 0$ and 
\[a(b(t,\cdot),\eta(t,\cdot)):=\int_M b(t,\cdot)\,\gamma(\eta(t,\cdot))\ dA\]
we let
 \begin{align*}
  (M_\eta b)(t,\cdot):= b(t,\cdot)- \psi\, \frac{a(b(t,\cdot),\eta(t,\cdot))}{a(\psi,\eta(t,\cdot))}.
 \end{align*}
\qed\smallskip

\begin{Lemma}\label{lemma:konvergenzen}
Let the sequence $(\eta_n)\subset Y^I$ satisfiy $\sup_n\norm{\eta_n}_{L^\infty(I\times
M)}<\alpha<\kappa$ and \eqref{eqn:schw}$_{(1,2)}$.
\begin{itemize}
 \item[\rm (1.a)] For $b\in
C(\bar I,L^2(M))$ the sequence $(\M_{\eta_n}b)$ converges to $\M_{\eta}b$ in $C(\bar I,L^2(M))$ independently of $\norm{b}_{C(\bar
I,L^2(M))}\le 1$.
\item[\rm (1.b)] Let $2\le r<\infty$. Provided that $b\in H^1(I,L^2(M))\cap
  L^r(I,H^2_0(M))$ and, additionally, that $(\pa_t\eta_n)$ converges in $L^2(I\times M)$ the sequence $(\M_{\eta_n}b)$ converges to $\M_\eta
  b$ in $H^1(I,L^2(M))\cap L^r(I,H^2_0(M))$.
\item[\rm (2.a)] Provided that $b\in C(\bar I,L^2(M))$ the sequence
  $(\F_{\eta_n}\M_{\eta_n}b)$ converges to $\F_{\eta}\M_{\eta}b$ in $C(\bar I,L^2(B_\alpha))$ independently of $\norm{b}_{C(\bar I,L^2(M))}\le 1$.
 \item[\rm (2.b)] On the conditions of $(1.b)$ the sequence $(\F_{\eta_n}\M_{\eta_n}b)$
converges to $\F_{\eta}\M_\eta b$ in
$H^1(I,L^2(B_\alpha))\cap L^r(I,W^{1,r}(B_\alpha))$.
\item[\rm (2.c)] Provided that $(b_n)$ converges to $b$ weakly in $L^2(I\times M)$ the sequence $(\F_{\eta_n}b_n)$ converges to
$\F_\eta b$ weakly in $L^2(I\times B_\alpha)$.
\end{itemize}
\end{Lemma}
\proof Comparing with \cite[Lemma A.11]{b101} only assertions $\rm (1.b)$ and $\rm (2.b)$ have changed. The proof of $\rm (1.b)$ proceeds exactly as before. The same is true for assertion $\rm (2.b)$ with the exception that, here, for the convergence in $L^r(I,W^{1,r}(B_\alpha))$ we have to exploit the fact that $Y^I$ embeds compactly into $L^\infty(I,W^{1,r}(M))$ which is a consequence of the classical Aubin-Lions lemma.
\qed\medskip

In the proof of Proposition \ref{lemma:komp} it comes in handy to have an $L^2$-orthogonal variant $\M_\eta^\perp$ of the operator $\M_\eta$ which is defined by
 \begin{align*}
  (\M_\eta^\perp b)(t,\cdot):= b(t,\cdot)- \gamma(\eta(t,\cdot))\, \frac{a(b(t,\cdot),\eta(t,\cdot))}{a(\eta(t,\cdot),\eta(t,\cdot))}.
 \end{align*}

\begin{Lemma}\label{lemma:mittelwertorth}
 For $\eta\in Y^I$ with $\norm{\eta}_{L^\infty(I\times M)}<\kappa$ the assertions of Lemma \ref{lemma:mittelwert} with $\M_\eta^\perp$ in place of $\M_\eta$ hold. Furthermore, for all $1\le r\le \infty$ and $0\le s\le 2$ we have 
\begin{equation*}
 \norm{\M_\eta^\perp b}_{L^r(I,H^s(M))}\le c\, \norm{b}_{L^r(I,H^s(M))}.
\end{equation*}
Finally, for $(\eta_n)\subset Y^I$ with $\sup_n\norm{\eta_n}_{L^\infty(I\times
M)}<\alpha<\kappa$ and \eqref{eqn:schw}$_{(1,2)}$ the claims {\rm (1.a)} and {\rm (2.a)} of Lemma \ref{lemma:konvergenzen} are true with $\M_\eta^\perp$ in place of $\M_\eta$.
\end{Lemma}
\proof The proofs proceed almost exactly as before. Note that $H^2(M)$ is an algebra. \qed\smallskip

\begin{Lemma}\label{lemma:ehrling} For  all $N\in\setN$, $6/5<p\le\infty$ and $\epsilon>0$ 
there exists a constant $c$ such that for all $\eta,\,\tilde\eta\in H^2_0(M)$ with
$\norm{\eta}_{H^2_0(M)}+\norm{\tilde\eta}_{H^2_0(M)}+\tau(\eta)+\tau(\tilde\eta)\le N$
and all $\bv\in W^{1,p}(\Omega_\eta)$, ${\tilde\bv\in W^{1,p}(\Omega_{\tilde\eta})}$ we have
\begin{equation*}
 \begin{aligned}
&\sup_{\norm{b}_{L^4(M)}\le
1}\bigg(\int_{\Omega_{\eta}} \bv\cdot\F_{\eta}\M_{\eta}b\ dx - \int_{\Omega_{\tilde\eta}}
\tilde\bv\cdot\F_{\tilde\eta}\M_{\tilde\eta}b\ dx\\
&\hspace{2cm} + \int_M \tr\bv\cdot\bnu\ \M_{\eta}b - \trt\tilde\bv\cdot\bnu\ \M_{\tilde\eta}b\
dA\bigg)\\
&\hspace{0.5cm}\le c\sup_{\norm{b}_{H^2_0(M)}\le
1}\bigg(\int_{\Omega_{\eta}} \bv\cdot\F_{\eta}\M_{\eta}b\ dx - \int_{\Omega_{\tilde\eta}}
\tilde\bv\cdot\F_{\tilde\eta}\M_{\tilde\eta}b\ dx\\
&\hspace{1cm}  + \int_M \tr\bv\cdot\bnu\ \M_{\eta}b -
\trt\tilde\bv\cdot\bnu\ \M_{\tilde\eta}b\ dA\bigg) + \epsilon\, \big(\norm{\bv}_{W^{1,p}(\Omega_{\eta})} +
\norm{\tilde\bv}_{W^{1,p}(\Omega_{\tilde\eta})}\big).
 \end{aligned}
\end{equation*}
\end{Lemma}
\proof The proof is a very simple modification of the proof of \cite[Lemma A.13]{b101} if we note that $W^{1-1/r,r}(M)$ embeds compactly into $L^{4/3}(M)$ for all $6/5<r<\infty$.
\qed
\medskip

\begin{Lemma}\label{lemma:projection}
Let $X$ be a function space that embeds compactly into $L^2(M)$, and let $\P_k$, $k\in\setN$, be the projection operators from the proof of Proposition \ref{lemma:komp}. Then for each $\eps>0$ we have
\[\norm{\P_k-id}_{\mathcal{L}(B,L^2(M))}\le\eps\]
provided that $k$ is sufficiently large.
\end{Lemma}
\proof By a simple compactness argument it suffices to show that for fixed $b\in L^2(M)$ we have
\[\norm{\P_k b-b}_{L^2(M))}\le\eps\]
provided that $k$ is sufficiently large. But this is an elementary consequence of the definition of the projection operators $\P_k$.\qed\smallskip

\section*{Acknowledgements}
The author would like to thank Lars Diening, Philipp N\"agele, and Michael R{\ocirc{u}}{\v{z}}i{\v{c}}ka for their valuable and helpful
comments and the fruitful discussions on the topic. The research of the
author was partly supported by the project C2 of the SFB/TR~71 ``Geometric Partial Differential Equations''.
%\end{acknowledgements}

% BibTeX users please use
%\bibliographystyle{myams}%spmpsci}
%\bibliography{biblio}   % name your BibTeX data base

\begin{thebibliography}{10}

\bibitem{b100}
{\sc G.~Acosta, R. Dur{\'a}n, and A. Lombardi}, {\sl Weighted {P}oincar\'e and {K}orn inequalities for {H}\"older {$\alpha$} domains}, Math. Methods Appl. Sci. {\bf 29} (2006), no. 4, 387--400.

\bibitem{b26}
{\sc G.~Acosta, R. Dur{\'a}n, and F. L{\'o}pez Garc{\'i}a}, {\sl Korn inequality and divergence operator: {C}ounterexamples and  optimality of weighted estimates}, Proc. Amer. Math. Soc. {\bf 141} (2013), no. 1, 217--232.

\bibitem{b37}
{\sc H.~Beirao~da Veiga}, {\sl On the Existence of Strong Solutions to a
  Coupled Fluid-Structure Evolution Problem}, Journal of Mathematical Fluid
  Mechanics {\bf 6} (2004), 21--52.

\bibitem{b38}
{\sc M.~Boulakia}, {\sl Existence of weak solutions for an interaction problem
  between an elastic structure and a compressible viscous fluid}, J. Math.
  Pures Appl. (9) {\bf 84} (2005), no.~11, 1515--1554.

\bibitem{b13}
{\sc F.~Boyer and P.~Fabrie}, {\sl {\'E}l{\'e}ments d'analyse pour l'{\'e}tude de
  quelques mod{\`e}les d'{\'e}coulements de fluides visqueux incompressibles},
  Math{\'e}matiques \& Applications (Berlin) [Mathematics \& Applications],
  vol.~52, Springer-Verlag, Berlin, 2006.

\bibitem{b109}
{\sc D. Breit, L. Diening, and S. Schwarzacher}, {\sl Solenoidal Lipschitz truncation for parabolic PDE's}, Preprint, arXiv:1209.6522v1 [math.AP].

\bibitem{b27}
{\sc A.~Chambolle, B.~Desjardins, M.J. Esteban, and C.~Grandmont}, {\sl
  Existence of weak solutions for the unsteady interaction of a viscous fluid
  with an elastic plate}, J. Math. Fluid Mech. {\bf 7} (2005), no.~3, 368--404.

\bibitem{b60}
{\sc C.H.A. Cheng, D.~Coutand, and S.~Shkoller}, {\sl Navier-{S}tokes equations
  interacting with a nonlinear elastic biofluid shell}, SIAM J. Math. Anal.
  {\bf 39} (2007), no.~3, 742--800 (electronic).

\bibitem{b42}
{\sc C.H.A. Cheng and S.~Shkoller}, {\sl The interaction of the 3{D}
  {N}avier-{S}tokes equations with a moving nonlinear {K}oiter elastic shell},
  SIAM J. Math. Anal. {\bf 42} (2010), no.~3, 1094--1155.

\bibitem{b44}
{\sc P.G. Ciarlet}, {\sl Mathematical elasticity. {V}ol. {II}}, Studies in
  Mathematics and its Applications, vol.~27, North-Holland Publishing Co.,
  Amsterdam, 1997, Theory of plates.

\bibitem{b45}
{\sc P.G. Ciarlet}, {\sl Mathematical elasticity. {V}ol. {III}}, Studies in Mathematics
  and its Applications, vol.~29, North-Holland Publishing Co., Amsterdam, 2000,
  Theory of shells.

\bibitem{b23} {\sc P.G. Ciarlet}, {\sl An introduction to differential
    geometry with applications to elasticity}, Springer, Dordrecht,
  2005, Reprinted from J. Elasticity {{\bf{7}}8/79} (2005), no. 1-3.

\bibitem{b40}
{\sc D.~Coutand and S.~Shkoller}, {\sl Motion of an elastic solid inside an
  incompressible viscous fluid}, Arch. Ration. Mech. Anal. {\bf 176} (2005),
  no.~1, 25--102.

\bibitem{b41}
{\sc D.~Coutand and S.~Shkoller}, {\sl The interaction between quasilinear elastodynamics and the
  {N}avier-{S}tokes equations}, Arch. Ration. Mech. Anal. {\bf 179} (2006),
  no.~3, 303--352.

\bibitem{b111}
{\sc G. Dal Maso and F. Murat}, {\sl Almost everywhere convergence of gradients of solutions to nonlinear elliptic systems}, Nonlinear Anal. {\bf 31} (1998), no. 3-4, 405--412.

\bibitem{b107}
{\sc L. Diening, J. M{\'a}lek, and M. Steinhauer}, {\sl On {L}ipschitz truncations of {S}obolev functions (with variable exponent) and their selected applications}, ESAIM Control Optim. Calc. Var. {\bf 14} (2008), no. 2, 211--232.

\bibitem{b108}
{\sc L. Diening, M. R{\ocirc{u}}{\v{z}}i{\v{c}}ka, and J. Wolf}, {\sl Existence of weak solutions for unsteady motions of  generalized {N}ewtonian fluids}, Ann. Sc. Norm. Super. Pisa Cl. Sci. (5) {\bf 9} (2010), no. 1, 1--46.

\bibitem{b106}
{\sc J. Frehse, J. M{\'a}lek, and  M. Steinhauer}, {\sl On analysis of steady flows of fluids with shear-dependent viscosity based on the {L}ipschitz truncation method},
SIAM J. Math. Anal. {\bf 34} (2003), no. 5, 1064--1083 (electronic).

\bibitem{b19}
{\sc G.P. Galdi, C.G. Simader, and H.~Sohr}, {\sl A class of solutions to
  stationary {S}tokes and {N}avier-{S}tokes equations with boundary data in
  {$W^{-1/q,q}$}}, Math. Ann. {\bf 331} (2005), no.~1, 41--74.

\bibitem{b28}
{\sc C.~Grandmont}, {\sl Existence of weak solutions for the unsteady
  interaction of a viscous fluid with an elastic plate}, SIAM J. Math. Anal.
  {\bf 40} (2008), no.~2, 716--737.

\bibitem{b110}
{\sc A. Hundertmark-Zau{\v{s}}kov{\'a}, M. Luk{\'a}{\v{c}}ov{\'a}-Medvi{\soft{d}}ov{\'a}, and G. Rusn{\'a}kov{\'a}}, {\sl Fluid-structure interaction for shear-dependent non-{N}ewtonian fluids}, Topics in mathematical modeling and analysis, Jind\u rich Ne\u cas Cent. Math. Model. Lect. Notes {\bf 7}, 109--158, Matfyzpress, Prague, 2012.

\bibitem{b54}
{\sc W.T. Koiter}, {\sl A consistent first approximation in the general theory
  of thin elastic shells}, Proc. {S}ympos. {T}hin {E}lastic {S}hells ({D}elft,
  1959) (Amsterdam), North-Holland, Amsterdam, 1960, pp.~12--33.

\bibitem{b20}
{\sc W.T. Koiter}, {\sl On the nonlinear theory of thin elastic shells. {I}, {II},
  {III}}, Nederl. Akad. Wetensch. Proc. Ser. B {\bf 69} (1966), 1--17, 18--32,
  33--54.

\bibitem{b102}
{\sc O. A. Lady{\v{z}}enskaja}, {\sl New equations for the description of the motions of viscous incompressible fluids, and global solvability for their boundary value problems},
Trudy Mat. Inst. Steklov. {\bf 102} (1967), 85--104.

\bibitem{b103}
{\sc  O. A. Lady{\v{z}}enskaja}, {\sl Modifications of the {N}avier-{S}tokes equations for large gradients of the velocities}, Zap. Nau\v cn. Sem. Leningrad. Otdel. Mat. Inst. Steklov. {\bf 7} (1968), 126--154.

\bibitem{b104}
{\sc O. A. Ladyzhenskaya}, {\sl The mathematical theory of viscous incompressible flow},  Second English edition, revised and enlarged. Mathematics and its Applications, Vol. 2, Gordon and Breach Science Publishers, New York, 1969.

\bibitem{b1}
{\sc J.M. Lee}, {\sl Introduction to smooth manifolds}, Graduate Texts in
  Mathematics, vol. 218, Springer-Verlag, New York, 2003.

\bibitem{b62}
{\sc D.~Lengeler}, {\sl Globale Existenz f{\"u}r die Interaktion eines
  Navier-Stokes-Fluids mit einer linear elastischen Schale}, Ph.D. thesis,
  Universit{\"a}t Freiburg, 2011, FREIDOK Server.

\bibitem{b101}
{\sc D.~Lengeler, M. R{\ocirc{u}}{\v{z}}i{\v{c}}ka}, {\sl Global
  weak solutions for an incompressible, Newtonian fluid interacting
  with a linearly elastic Koiter shell}, Preprint, arXiv:1207.3696v1 [math.AP].

\bibitem{b61}
{\sc J.~Lequeurre}, {\sl Existence of strong solutions to a fluid-structure
  system}, SIAM J. Math. Anal. {\bf 43} (2011), no.~1, 389--410. 

\bibitem{b105}
{\sc J.-L. Lions}, {\sl Quelques m\'ethodes de r\'esolution des probl\`emes aux
limites non lin\'eaires}, Dunod, 1969.

\end{thebibliography}

\ifx\undefined\bysame
\newcommand{\bysame}{\leavevmode\hbox to3em{\hrulefill}\,}
\fi

% \begin{thebibliography}{3}
% %
% % and use \bibitem to create references. Consult the Instructions
% % for authors for reference list style.
% %
% % Format for Journal Reference
% \bibitem{Ref1}
% Author, I.: Article title. Journal Title-Abbreviated {\bf Vol}, pp--pp (year)
% \bibitem{Ref2}
% Author, I., Smith, J.: Book Title. Publisher, Place (year)
% % Format for proceedings
% \bibitem{Ref3}
% Author, I., Smith, J.: Paper title. In: Editor, A. (ed.) Proceedings
% Title, Location, Date, pages. Publisher, Place (year)
% \end{thebibliography}

\end{document}